\documentclass[preprint,12pt]{elsarticle}
\usepackage[T1]{fontenc}

\usepackage{amssymb}
 \usepackage{amsthm}
 \usepackage{amsmath}
 \usepackage{tikz}
 \usepackage{caption,subcaption}
 \usepackage{hyperref}
 \usepackage{comment}
 \usepackage{enumitem}

 \newtheorem{theorem}{Theorem}

\newtheorem{corollary}[theorem]{Corollary}

\newtheorem{conjecture}[theorem]{Conjecture}

\theoremstyle{definition}

\newcommand{\TT}{\mathbb{T}}
\newcommand{\NN}{\mathbb{N}}
\newcommand{\ZZ}{\mathbb{Z}}
\newcommand{\QQ}{\mathbb{Q}}
\newcommand{\RR}{\mathbb{R}}
\newcommand{\FF}{\mathbb{F}}

\newcommand{\Rab}{\textrm{Rab}}

\newcommand{\LR}{\textrm{LR}}
\newcommand{\LRC}{\textrm{LR Conjecture}}
\newcommand{\LRP}{\textrm{LR Problem}}

\journal{Computer Science Review}

\begin{document}

\begin{frontmatter}



\title{The Lonely Runner Conjecture turns 60\tnoteref{label3}}
\tnotetext[label3]{The authors were supported by the Spanish Agencia Estatal de Investigación under projects PID2020-113082GB-I00, PID2023-147202NB-I00, RED2022-134947-T and the Severo Ochoa and María de Maeztu Program for Centers and Units of Excellence in R\&{}D (CEX2020-001084-M).}


\author[label1,label2]{Guillem Perarnau}
\author[label1,label2]{Oriol Serra}

\affiliation[label1]{organization={Department of Mathematics, Universitat Politècnica de Catalunya},
city={Barcelona},
            country={Spain}
            }
\affiliation[label2]{organization={Centre de Recerca Matemàtica},
city={Barcelona},
            country={Spain}
            }
%
%

\begin{abstract}
The Lonely Runner Conjecture originated in Diophantine approximation is turning 60. Even if the conjecture is still widely open, the flow of partial results, innovative tools and connections to different problems and applications has been steady on its long life. This survey attempts to give a panoramic view of the status of the problem, trying to highlight the contributions of the many papers that it has originated.
\end{abstract}

%

\begin{keyword} Lonely Runner Conjecture \sep Diophantine approximation \sep View-Obstruction problems \sep Chromatic number of distance graphs



\end{keyword}

\end{frontmatter}

\tableofcontents



\section{Introduction}
\label{sec:intro}

The \emph{Lonely Runner (\LR) Conjecture} was posed  by Wills~\cite{W1968} in 1968 (see also~\cite{W1965,W1967}), in the context of Diophantine approximation,  and by
Cusick~\cite{C1973} in 1973, in the context of view--obstruction problems. Its picturesque name comes from the following interpretation due to Goddyn~\cite{BGGST1998}.
Consider a set of {$n$} runners  on the unit circle running with different constant speeds and
starting at the origin. The conjecture states that, for each runner, there is a time where she is \emph{lonely}, precisely, at distance at least {$1/n$} on the circle from all the other
runners.

More precisely, the conjecture can be stated as follows. For any real number $x$, denote by $\|x\|$ the distance from $x$ to the closest
integer
$$
\|x\| = \min \{x-\lfloor x\rfloor,  \lceil x\rceil -x\}.
$$

{
\begin{conjecture}[\LRC, Original version]\label{conj:LRC_original}
For every $n\in \NN$, every $i\in [n]$ and every set of pairwise distinct  real numbers $v_1,\dots ,v_{n}$, there exists $t\in
\RR$ such that
\begin{equation}\label{eq:ALSD}
\min_{j\neq i} \|t(v_j -v_i)\|\geq \frac{1}{n}.
\end{equation}
\end{conjecture}
}

Since Conjecture~\ref{conj:LRC_original} is invariant by translations of the speed set $\{v_1,\dots,v_n\}$, it is convenient to change the point of reference so one of the runners has zero speed. This runner is then omitted from the problem and the conjecture is formulated in the following simpler form.

\begin{conjecture}[\LRC]\label{conj:LRC}
For every $n\in \NN$ and every set of nonzero real numbers $v_1,\dots ,v_{n}$, there exists $t\in
\RR$ such that
\begin{equation}\label{eq:SLDO}
 \min_{1\leq i \leq n}\| tv_i\|\geq \frac{1}{n+1}.
\end{equation}
\end{conjecture}
{Since the stagnant runner at the origin is omitted, proving this conjecture for a given $n\in \NN$ proves the original one for $n+1$ runners. Observe that the second version does not need the speeds to be distinct, but we require them to be nonzero. In contrast to the original one, Conjecture~\ref{conj:LRC} is not invariant by translations but it is invariant  by dilations.}


The purpose of this paper is to survey the state-of-the-art on the \LRP. To the best of our knowledge, this is the first in-depth review of the problem, see \cite{SW2018} for a shorter overview intended for a broader audience. After introducing some reductions in Section \ref{sec:reductions} that will be adopted throughout the paper, we begin in Section \ref{sec:connections} by reviewing various contexts in the literature where the problem has emerged --- from Diophantine approximation to Geometry, Graph Theory or Arithmetics --- which provides a broader perspective on the problem. Section \ref{sec:tight} addresses the tightness of the conjecture by exploring sets of speeds for which the equality holds. Section~\ref{sec:gap} reviews the efforts to increase the so-called \emph{gap of loneliness}, from trivial bounds to more sophisticated arguments. There has been some instructive progress on settling the conjecture for small values of $n$, which is discussed in Section~\ref{sec:small}. The conjecture has also been confirmed for several classes of speed sets and these are detailed in Section~\ref{sec:special}. Section~\ref{sec:bounded} is devoted to a result by Tao, which reduces the problem to speeds bounded by a function of $n$ and, conversely, it also shows progress when all speeds are linear in $n$. Random sets of runners satisfy a much stronger inequality than the one predicted in Conjecture~\ref{conj:LRC}, the topic is discussed in Section~\ref{sec:random}. In Section~\ref{sec:variations} we explore several extensions of the problem, including invisible runners and the shifted conjecture, among others. The paper concludes in Section~\ref{sec:final} with final remarks.

\section{Reductions}\label{sec:reductions}

The statement of the \LRC\ as in Conjecture~\ref{conj:LRC} can be simplified by a number of observations and auxiliary results. First of all, observe that we may assume that all speeds are positive, as otherwise we may replace $v_i$ by $-v_i$, and the set of times for which~\eqref{eq:SLDO} is satisfied does not change.  For the sake of simplicity, we will assume that the speeds are ordered increasingly
$$
0<v_1<v_2<\dots< v_n,
$$
and denote $V=\{v_1,\dots, v_n\}$ for the set of speeds.
Eventually, and particularly in the geometric interpretations of the problem, we will drop the increasing speeds assumption and consider them as a vector instead of a set, denoted by $\mathbf{v}=(v_1,v_2,\dots, v_n)$.

Originally, Wills~\cite{W1968} proposed a version of the \LRC\, with integral speeds (see Section~\ref{subsec:diophantine}) and Cusick~\cite{C1973} claimed 
that it was equivalent to Conjecture~\ref{conj:LRC}.
For that, one can appeal to Kronecker's classical theorem on Diophantine approximation (see e.g.~\cite{HW1979}), from which we state a generalisation.
\begin{theorem}[Kronecker's Theorem]\label{thm:kron}
For every $n\in \NN$, every $\epsilon>0$, every $a_1,a_2,\dots, a_n \in \RR$ and every $ v_1, v_2, \dots, v_n\in \RR$ that are linearly independent over $\QQ$, 
there exists $t\in \RR$ such that
\begin{equation}\label{eq:EUDS}
\max_{1\leq i\leq n}\| tv_i-a_i\| \leq \epsilon \quad \text{for every $i\in [n]$}\;.    
\end{equation}
In particular, if $v_1,v_2,\dots, v_n$ have rational dimension $m\leq n$, there exists an $m\times n$ matrix $A$ with rational entries of rank $m$ such that for each $\mathbf{a}=(a_1,a_2,\dots,a_n)\in \text{Ker}(A)+\ZZ^n$ there exists $t\in \RR$ such that~\eqref{eq:EUDS} holds.
\end{theorem}
If $V=\{v_1,v_2, \dots, v_n\}$ is the set of speeds and $a_i=1/2$ for $i\in [n]$, then the first part of the theorem ensures the existence of a time when all the runners are arbitrarily far from the origin, provided that $V$ is linearly independent over $\QQ$. This conclusion is much stronger than the desired one and it was used by Bohman, Holzman and Kleitman~\cite[Section 4]{BHK2001} to reduce the problem to rational speeds (and thus to integral ones) as follows: if all speeds are not collinear over $\QQ$, then $\dim (\text{Ker} (A)) =m\geq 2$ and we can find two linearly independent vectors with rational entries. These can be used to construct another vector with rational entries where two of the entries are equal but of opposite sign. This reduces the problem to $n-1$ speeds with rational entries and the same gap.
Recently, Henze and Malikiosis~\cite[Section 5]{HM2017} provided an alternative proof of the reduction that does not rely on the lower dimensions, see also~\cite{DT2020}.

Once the problem has been reduced to positive integer speeds, it suffices to consider $t$ taking values in $(0,1)$. Indeed, $\|tv_i\|=0$ for any $i\in [n]$ and $t\in \ZZ$, so the vector of runner's positions as a function of $t$ is periodic of period at most $1$. 

\begin{conjecture}[\LRC; Integer version]\label{conj:LRC_int}
For every $n\in \NN$ and every $n$-set $V$ of positive integers, there exists $t\in
(0,1)$ such that 
\begin{equation*}
\min_{v\in V} \| tv\|\geq \frac{1}{n+1}\;.
\end{equation*}
\end{conjecture}
If needed, one may restrict the problem to speed sets for which $\gcd(V)=1$, as the problem is invariant by dilations. For further reductions of the problem, see Section~\ref{sec:simple_rem}.

Finally, one can ensure that, if the origin becomes isolated, then it does at certain times. Namely, if $t_0\in (0,1)$ is a time $t$ for which the supremum of $\min_{v\in V}\|tv\|$ is attained, then, by continuity, that time must minimize simultaneously the distance to the origin of two of the runners. Suppose $v_i, v_j$ are the speeds of such runners. Then $t_0v_i=1-t_0v_j \mod 1$ and we can write
\begin{equation}\label{eq:time_red}
t_0=\frac{\ell}{v_i+v_j},
\end{equation} 
for some $\ell\in \NN$. This is an important observation which has been repeatedly used in the literature.

\section{Connections}\label{sec:connections}

Several connections of the \LRP\ with different areas have been established in the literature. This wide range of connections not only enriches the view on the problem but also  provides  tools and techniques to get new results. In the rest of the survey, we will often consider these different interpretations of the problem so it is appropriate to discuss them here.

\subsection{Diophantine approximation}\label{subsec:diophantine}  

The \LRP\ first arose in connection to Diophantine approximation. For any $V\subset \RR$, let
\begin{equation}\label{eq:xi_kappaV}
\begin{aligned}
\xi (V) & =\sup_{t\in \ZZ}\min_{v\in V}\|tv\|,\\
\kappa (V) & =\sup_{t\in (0,1)}\min_{v\in V}\|tv\|.
\end{aligned}
\end{equation}
For any positive integer $n$, Wills~\cite{W1968} defines
\begin{align*}
\xi(n)&=\inf_{V\subset \RR\setminus Q \atop |V|=n}\xi(V),\\
\kappa (n)&=\inf_{V\subset \NN \atop |V|=n} \kappa(V),
\end{align*}
as extensions of the Diophantine approximation parameters, for a  simultaneous homogeneous Diophantine approximation problem \footnote{In fact,  Wills uses $\kappa$ to denote what we call by $\xi$; however, the notation used here has become the standard in the literature.}; see also \cite{W1967}, where Wills studied $\xi(n)$ without explicitly define it. Indeed, traces of the problem can be found as early as Wills' 1965 doctoral thesis~\cite{W1965}, which inspired the title of this survey. We informally refer to $\kappa(V)$ and $\kappa(n)$ as the \emph{gaps of loneliness}.
Wills~\cite{W1968} showed that $\kappa(n)=\xi(n)$.

Using this terminology, the \LRC\ can be phrased as
$$
\inf_{V\subset \RR\setminus\{0\} \atop |V|=n} \kappa(V)\geq \frac{1}{n+1} \quad\text{for all }n\in \NN.
$$
Since the \LRP\ with real speeds can be reduced to the \LRP\ with integral speeds, as shown in the previous section, we obtain the following compact formulation of the \LRC.
\begin{conjecture}[\LRC; Diophantine interpretation]\label{conj:LRC_dioph}
For every $n\in \NN$, we have
\begin{equation*}
\kappa(n)\geq \frac{1}{n+1}\;.
\end{equation*}
\end{conjecture}


Wills~\cite{W1968} established the following inequalities 
\begin{equation}\label{eq:SODP}
\frac{1}{2n}\le \kappa (n)\le \frac{1}{n+1}.
\end{equation}
{See Sections~\ref{sec:tight} and~\ref{sec:gap} respectively for a detailed discussion on the upper and on the lower bound.}

One can also consider the following discrete version of the problem. For a positive integer $N$ and $v\in \ZZ/N\ZZ$, let us define $\|v\|_N$ as the circular distance from $v$ to zero in $\ZZ/N\ZZ$. This is the residue class of  $v$ or of $-v$ modulo $N$ in $\{0,1,\ldots ,\lfloor N/2\rfloor\}$. For a set $V$ of positive integers, we abuse notation and still denote by $V$ the set of residue classes of the elements of $V$ modulo $N$. Define
\begin{equation}\label{def:kappa2}
\kappa_N(V)=\max_{\lambda\in (\ZZ/N\ZZ)^*}\min_{v\in V}\|\lambda v\|_N
\end{equation}
{where $(\ZZ/N\ZZ)^*$ denotes the multiplicative group of the integers modulo $N$.} If there is $\lambda \in (\ZZ/N\ZZ)^{{*}}$ such that
$$
\|\lambda v\|_N\ge N/k\quad \text{for all}\; v\in V,
$$
then $\|(\lambda/N)v\|\ge 1/k$. Therefore, 
$$
\kappa (V)=\limsup_{N\to\infty} \frac{1}{N}\kappa_N(V).
$$
As remarked by Haralambis~\cite[Remark 1]{H1976b} (see also Czerwi\'{n}ski and Grytczuk~\cite[Theorem 6]{CG2008} for a proof), {using the observation made in~\eqref{eq:time_red}}, the $\limsup$ can be replaced by the following finite optimization
\begin{equation}\label{eq:SPDO}
\kappa (V)=\max_{N=v+v'\atop v,v'\in V} \frac{1}{N}\kappa_N(V).
\end{equation}
In particular, for a given $n$-set $V$, the \LRC\ can be verified to hold for $V$ in $O(n^2v_n)$ time\footnote{The fact that $\kappa(V)$ is computable in polynomial time for every given $V$ clearly follows from the fact that $\sup_{t\in(0,1)} \min_{v\in V}  \|tv\|$ is a piecewise linear function in $t$.}.
 
\subsection{View--Obstruction and Billiard Trajectory problems} \label{subsec:VO_Bill}

Shortly after Wills' paper, Cusick~\cite{C1973} introduced the following geometric problem, which he named the \emph{View--Obstruction Problem}. Let $K$ be a convex body in $\RR^n$ containing the origin. Denote by $\lambda K$ the dilation of $K$ by $\lambda\in \RR^+$. Let $\mathbf{1}\in \RR^n$ be the all ones vector. Let 
$$
\Delta (K,\lambda)=\lambda K+ \NN^n - \tfrac{1}{2}\mathbf{1},
$$ 
{be the set of all translates of $\lambda K$ centered at points with positive half-integer coordinates}. Then the problem is to estimate $\lambda_K$, the minimum value of $\lambda$ such that, for every vector $\mathbf{v}=(v_1,\ldots ,v_n)$ of positive real numbers, the ray $R_\mathbf{v}=\{(x_1,\ldots ,x_n)=(tv_1,\ldots ,tv_n):\, t\in \RR^+\}$ intersects $\Delta(K,\lambda)$. 

If $K=Q_n=[0,1]^n$ is the unit cube, we denote $\lambda (n)=\lambda_{Q_n}$. As shown in~\cite[Lemma 1]{C1973}, we have 
$$
\lambda (n)=1-2\kappa (n).
$$
Hence, the above Diophantine approximation problem is equivalent to the View-Obstruction problem for the cube. 
{
\begin{conjecture}[\LRC, View-Obstruction interpretation]\label{conj:LRC_VO}
For every $n\geq 1$, we have
\begin{equation}\label{eq:PSKD}
\min \left\{\lambda >0:\, \{\mathbf{v}:\, R_\mathbf{v}\cap \Delta(Q_n,\lambda)\neq \emptyset \}\supset \NN^n\right\} =\frac{n-1}{n+1}.
\end{equation}
\end{conjecture}
}

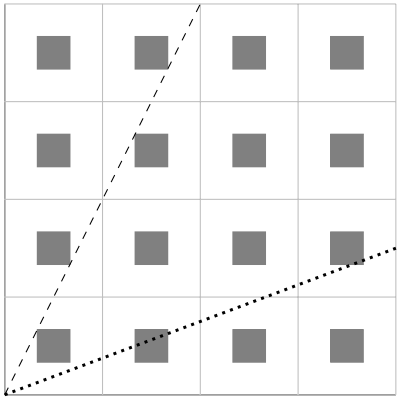
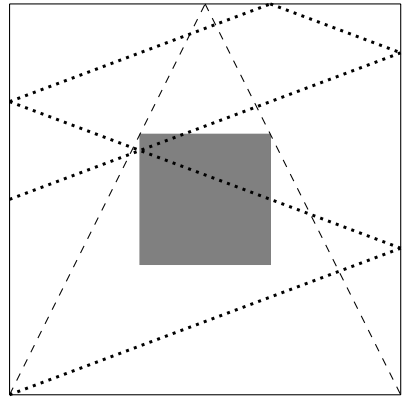
\begin{figure}[h]

\begin{subfigure}[b]{0.4\textwidth}
\begin{center}
\begin{tikzpicture}[scale=1.3]
\draw (0,0)--(4,0) (0,0)--(0,4);
\draw[lightgray] (0,0) grid (4,4);
\foreach \i in {1,...,4}
{
\foreach \j in {1,...,4}
{
\draw[fill,gray] (\i-1/2-1/6,\j-1/2-1/6) rectangle  (\i-1/2+1/6,\j-1/2+1/6);
}
}
\draw[line width=0.4mm,dotted] (0,0)--(4,3/2);
\draw[dashed] (0,0)--(2,4);
\end{tikzpicture}
\end{center}
\caption{View-obstruction interpretation: all rays eventually intersect the translated dilated squares.}  
    \end{subfigure}
\begin{subfigure}[b]{0.15\textwidth}
\hspace{1cm}
\end{subfigure}
\begin{subfigure}[b]{0.4\textwidth}
      \begin{center}
\begin{tikzpicture}[scale=1.3]
\draw (0,0)--(4,0) (0,0)--(0,4) (0,4)--(4,4) (4,0)--(4,4);
\draw[fill,gray] (2-2/3,2-2/3) rectangle  (2+2/3,2+2/3);
\draw[dashed] (0,0)--(2,4);
\draw[dashed] (2,4)--(4,0);
\draw[line width=0.4mm,dotted] (0,0)--(4,3/2);
\draw[line width=0.4mm,dotted] (4,3/2)--(0,3);
\draw[line width=0.4mm,dotted] (0,3)--(4*2/3,4);
\draw[line width=0.4mm,dotted] (4*2/3,4)--(4,4-1/2);
\draw[line width=0.4mm,dotted] (4,4-1/2)--(0,2);
\end{tikzpicture}
\end{center}
      \caption{Billiard ball trajectories interpretation: all ball trajectories eventually step on the shaded square.}
\end{subfigure}
\caption{Illustration for $n=2$ and $\lambda=1/3$. The dashed line corresponds to a case where equality is attained.}\label{fig:view}
\end{figure}

A simple geometric argument solves the case $n=2$, with $\lambda (2)= 1/3$, see an illustration in Figure~\ref{fig:view}. An argument is given in~\cite{C1973} for the case $n=3$. 

In this form, the problem can be also interpreted as a \emph{Billiard Ball Trajectory} problem, see Schoenberg~\cite{S1976}. Consider $Q_n$ and paint a smaller cubic volume of side length $\lambda\in (0,1)$ centered at $(1/2)\mathbf{1}$. Place a ball at the origin and shoot it in any positive direction. When the ball hits one of the walls of the cube, it bounces off the wall with perfect symmetry. Then, $\lambda(n)$ is the minimum of the values of $\lambda$ for which the trajectory of the ball eventually steps on the painted cube regardless of the initial direction. For $n=2$, see an illustration in Figure~\ref{fig:view}.

This perspective is also suited to discuss a natural variant of the
problem where the ball can be initially placed at any point of $Q_n$, or in the original jargon, where the runners start at different points on the track. This is known as the shifted variant of the \LRC; see Section~\ref{subsec:shifted}.

\subsection{Covering radius of zonotopes}\label{subsec:zonotope}

An equivalent geometric formulation of the \LRC\ in terms of covering radius of polytopes, a classical topic in Discrete Geometry and Geometry of Numbers, is discussed by Henze and Malikiosis~\cite{HM2017}. 

A \emph{lattice zonotope} $Z$ in $\RR^m$ is the Minkowski sum of $n\ge m$ segments of the form {$[0,\mathbf{z}_1]+\cdots +[0,\mathbf{z}_n]$, where $\mathbf{z}_1,\ldots ,\mathbf{z}_n\in \ZZ^m$ are the \emph{generators} of the polytope. The \emph{center} of $Z$ is defined as $\mathbf{x}=\tfrac{1}{2}\sum_{i=1}^n \mathbf{z}_i$}. We will assume that the points {$\mathbf{z}_1,\ldots ,\mathbf{z}_n$} are in general position, namely, every subset of $m$ points forms a basis of $\RR^m$. We fix $m=n-1$.

With this terminology the following conjecture, closely related to the View-Obstruction problem, is shown in~\cite{HM2017} to be equivalent to the \LRC.
\begin{conjecture}[\LRC; Polytopal interpretation 1] For every $n\in \NN$ and every lattice zonotope $Z$ generated by $n$ vectors in general position in $\ZZ^{n-1}$, we have
$$
\frac{n-1}{n+1}(Z-\mathbf{x})\cap \ZZ^{n-1}\neq \emptyset,
$$
where $\mathbf{x}$ is the center of the zonotope.
\end{conjecture}
This geometric  approach also provides a direct way to prove that the \LRC\ can be reduced to integer speeds; see~\cite[Lemma 5.3]{HM2017}. 

A closely related geometric interpretation is developed by Beck and Schymura~\cite{BS2024} exploiting the symmetry properties of a body. The \emph{coefficient of asymmetry} of a convex body $K$ with respect to an interior point $\mathbf{w}$, denoted by $\alpha_{K,\mathbf{w}}$, is defined as the minimum value of $\alpha\ge 1$ such that $\mathbf{w}-K$ is contained in the dilation $\alpha(K-\mathbf{w})$. For instance, $\alpha_{K,\mathbf{w}}=1$ if and only if $K$ is centrally symmetric around $\mathbf{w}$. The \LRC\ is equivalent to the following one.

Given $\mathbf{v}\in \ZZ^n$, define the zonohedron
$$
Z_{\mathbf{v}}= \mathbf{v}\RR + [0,1]^n.
$$

\begin{conjecture}[\LRC; Polytopal interpretation 2] For every $n\in \NN$ and every $\mathbf{v}\in \NN^n$ there is an interior lattice point $\mathbf{w}\in Z_{\mathbf{v}}\cap \ZZ^n$ such that 
$$
\alpha_{Z_{\mathbf{v}},\mathbf{w}}\leq n.
$$
\end{conjecture}

The geometric interpretation in \cite{BHS2019} is in between the view obstruction and the zonotopes. For an integral vector $\mathbf{m}\in \ZZ^k_+$ the {\it polyhedral cone} $K_{\mathbf{m}}$ is the cone spanned by the vectors of the form $(k+1)\mathbf{m}+\mathbf{u}$ for all $\mathbf{u}\in \{1,k\}^k$. It is shown in \cite{BHS2019} that, for a vector $\mathbf{v}=(v_1,\ldots ,v_k)$ of velocities,   the polyhedron $P_{\mathbf{v}}=\mathbf{v} \mathbb{R}+[-\frac{1}{k+1},\frac{1}{k+1}]^k$ intersects the lattice (the view obstruction approach) if and only if there is $\mathbf{m}\in \ZZ^k_+$ such that $\mathbf{v}\in K_{\mathbf{m}}$. The additional value of this polyhedral interpretation is that the vector $\mathbf{m}$ is connected to the time when the runners get lonely: if  $\mathbf{v}\in K_{\mathbf{m}}$ then the $i$--th runner is in its $(m_i+1)$--th lap around the track at the time when they are all far from the origin.

One feature of the above geometric interpretations of the \LRC\ is that they can be adapted to describe a generalization of the conjecture where the runners can start at different positions on the circle; see Subsection \ref{subsec:shifted}.

\subsection{Nowhere zero flows} Bienia, Goddyn, Gvozdjak, Seb\H o,  and Tarsi~\cite{BGGST1998} established a further connection between the \LRP\ and flows in graphs. A \emph{nowhere zero flow} in an undirected graph $G$ is an orientation of the edges of $G$ together with a function $f$ from the edge set to $\NN$ such that, at each vertex $v$, the sum of the values of $f$ on the edges entering $v$ equals the sum of the values on edges leaving $v$. The main result in~\cite{BGGST1998} is the following one.


\begin{theorem}\label{thm:flow} If a graph $G$ admits a nowhere zero flow with $n$ distinct values then it also admits a nowhere zero flow with values in $[n]$.
\end{theorem}

{The authors set the following connection between this result and the \LRP. Suppose that $f$ is a nowhere zero flow in a graph $G$ taking $n$ distinct values $f_1,\ldots ,f_n$. If the \LRC\ holds, then there exists $t\in (0,1)$ such that $tf$ is a nowhere zero flow and $\|tf_i\|\geq \frac{1}{n+1}$ for each $i\in [n]$. By a theorem of Ford and Fulkerson~\cite{FF1958}, there exists an integer flow, denoted by $\lfloor tf\rceil$, which takes values $\lfloor tf\rfloor$ or $\lceil tf\rceil$ on the edges of $G$. Thus, $f'= (n+1)(f-\lfloor tf\rceil)$ is a nowhere zero flow taking values in $[-n,-1]\cup[1,n]$, Again, by~\cite{FF1958}, there exists a nowhere zero flow taking values in $[n]$.} 

The proof of Theorem~\ref{thm:flow} relies on the validity of the \LRC\ for $n\le 4$, for which the authors give a simpler proof for the case $n=4$, and on the celebrated six--flow theorem by Seymour~\cite{S1981} for $n\ge 5$. Nevertheless, by the above argument, proving the \LRC\ for each $n$ would also imply the validity of Theorem~\ref{thm:flow}.

The statement of Theorem~\ref{thm:flow} can be generalised to the setting of regular matroids, which can be defined through totally unimodular matrices. Recall that a matrix $A$ is \emph{totally unimodular} if every subdeterminant belongs to $\{-1,0,1\}$. The following open problem can be understood as a weak version of the \LRC:


\begin{conjecture}[\cite{BGGST1998}]\label{conj:matroid} Let $A$ be a totally unimodular $k\times m$ matrix. If there is a solution $\mathbf{x}=(x_1,\ldots ,x_m)$ of the equation $A\mathbf{x}=0$ with nonzero entries and {$|\{x_i:\, i\in [m]\}|\leq n$}, then there is also a solution $\mathbf{x}'=(x'_1,\ldots ,x'_m)$ with $|x'_i|\in [n]$ for all $i\in [m]$.
\end{conjecture}

{The conjecture implies Theorem~\ref{thm:flow} by taking $A$ to be the vertex-edge incidence matrix of $G$.}

\subsection{Chromatic number}\label{subsec:chrom}

Motivated by the plane coloring problem, the minimum number of colors needed to color the points in $\RR^2$ such that points at distance one receive distinct colors, Eggleton, Erd\H{o}s and Skilton~\cite{EES1985} introduced the study of the chromatic number, {$\chi(G)$}, of distance graphs. A particularly interesting case are the graphs on the integers, $G(\ZZ,D)$, where two vertices are joined if and only if their difference in absolute value belongs to a prescribed set $D=\{d_1,\ldots ,d_n\}$ of distances. 

{Since $-n$ and $n$ have degree at most $|D|$ in the subgraph induced by $[-n,n]$ in $G(\ZZ,D)$, one can color the graph greedily with $|D|+1$ colours and
\begin{equation}\label{eq:DPSI}
\chi(G(\ZZ,D))\leq |D|+1.
\end{equation}
}

In order to place the \LRC\;  in this setting,  two related parameters arise, the fractional and the circular chromatic numbers, which provide an approach to the study of the chromatic number of distance graphs. Let $\mathcal{I}(G)$ denote the family of independent sets of a graph $G$ and, for a vertex $x$ of $G$, let ${\mathcal I}(G,x)$ denote the family of independent sets containing $x$. The \emph{chromatic number} of a graph $G$ can be defined as the minimum {non-negative} integer solution of $\sum_{I\in {\mathcal I}(G)}x_I$ subject to $\sum_{I\in {\mathcal I}(G,x)}x_I\ge 1$ for each vertex $x$ of $G$. The relaxation of this Integer Linear Program gives rise to the \emph{fractional chromatic number} $\chi_f(G)$ of the graph. By definition, $\chi_f(G)$ provides a lower bound to $\chi(G)$, and it is lower bounded by $\chi(G)-1$.

The classical \emph{Motzkin problem} in combinatorial number theory asks for the maximum asymptotic density $\mu(D)$ of a set $S\subset \NN$ which avoids a prescribed set $D=\{d_1,\ldots ,d_n\}$ of differences, namely, $|x-y|\not\in D$ for every $x,y\in S$. Chang, Liu and Zhu~\cite{CLZ1999} proved that $\mu(D)$ is equivalent to the fractional chromatic number of the distance graph $G(\ZZ,D)$:
$$
\chi_f(G(\ZZ,D))=\frac{1}{\mu (D)}.
$$
The \emph{circular chromatic number} $\chi_c(G)$ of a graph $G$ is the minimum value of $r\in \RR^+$ such that there is a map $f$ from the vertex set to $(0,1)$ such that, for every edge $xy$ of $G$, the distance between $f(x)$ and $f(y)$ in the unit torus $\RR/\ZZ$ is at least $1/r$. It turns out that 
\begin{equation}\label{eq:FDID}
\chi_f(G)\le \chi_c(G)\le \chi (G),
\end{equation}
for any locally finite countable graph $G$, see e.g. Zhu~\cite{Z2001} for a thorough discussion on this parameter. The definition given above is reminiscent of the \LRP, and in fact it is shown there that
$$
\chi_c(G(\ZZ,D))\le \frac{1}{\kappa (D)},
$$
{where $\kappa(D)$ is defined as in~\eqref{def:kappa2}.}
This gives the chain of inequalities
$$
\chi (G(\ZZ, D))-1 \leq \frac{1}{\mu (D)}=\chi_f(G(\ZZ, D))\le \chi_c(G(\ZZ, D))\le \frac{1}{\kappa (D)},
$$ 
connecting the study of the chromatic numbers of distance graphs with Motzkin problem and the \LRP, see e.g.  Liu~\cite{L2008}. {From~\eqref{eq:DPSI} and~\eqref{eq:FDID}, one deduces that $\chi_c(G)\leq |D|+1$, while the \LRC\ asks whether $1/\kappa(D)\leq |D|+1$. This can be seen as additional support to the conjecture.}
The inequality $\kappa (D)\le \mu (D)$ is often strict but there are cases of equality. {In particular this has been proved for $|D|=2$ (see Section~\ref{sec:small}), is conjectured for $|D|=3$~\cite{L2008} and false for $|D|=4$~\cite{LZ2004}.} 

The connection with the circular chromatic number is also linked with the  discrete version of the \LRP\ discussed at the end of Section~\ref{subsec:diophantine}. From~\eqref{eq:SPDO}, we can write
$$
\kappa (D)=\max_{N=d+d'\atop d,d'\in D} \frac{1}{N}\kappa_N(D).
$$
In the graph language, for any such $N$, the distance graph $G(\ZZ,D)$ admits a graph homomorphism onto the circulant graph $G(\ZZ/N\ZZ, D)$. Therefore, any chromatic number of the latter is a lower bound for the chromatic number of the former, see Liu~\cite{L2008} for more details on this perspective.

\subsection{Coprime mappings}\label{subsec:coprime_map}

A more recent connection targeted to sets of small speeds is established by Bohman and Peng~\cite{BP2022} via coprime mappings. Given two finite sets $A,B$ of integers with the same cardinality, a bijection $f:A\to B$ is a \emph{coprime mapping} if $a$ and $f(a)$ are coprime for every $a\in A$. Pomerance and Selfridge~\cite{PS1980} proved that there is a coprime {mapping} from $[n]$ to any set $B$ of $n$ consecutive integers, solving a conjecture of Newman. In~\cite[Theorem 2.1]{BP2022}  the authors prove a weaker version which, for the purpose of its application to the \LRC, can be stated as follows.

\begin{theorem}\label{thm:cop} Let $I,J\subset [n]$ be two sets of $2m$ consecutive integers. There is an absolute constant $c>0$ such that, if  $m\ge e^{c(\log\log n)^2}$, then, for every nonempty subsets $S\subset I$ and $ T\subset J$ with $|S|+|T|\ge 2m$ not consisting entirely of even numbers, there are $s\in S$ and $t\in T$ which are coprime.
\end{theorem}

In particular, by Hall's theorem, if $S$ and $T$ satisfy the conditions of Theorem~\ref{thm:cop} and in addition $|S|+|T|\ge 2m+1$ then there is a coprime mapping from $I$ to $J$. 

The connection with the \LRP\ is as follows. Let $s,t$ be two coprime numbers and $m=s+t$, so that $s\in (\ZZ/m\ZZ)^*$ and $q=s^{-1}$. Then, $q/m$ is an appropriate time for every set of speeds which are not in the dilated interval $(q/m)\cdot[\tfrac{1}{n+1},\tfrac{n}{n+1}]$. This argument is enough to solve the \LRC\ when all speeds are small; see Section~\ref{sec:bounded}.



\section{Tight Instances}\label{sec:tight}

Qualitatively speaking, the \LRC\ states that, eventually, all the runners are far away from the origin, leaving it lonely. From the quantitative point of view, it ventures that a gap of $\tfrac{1}{n+1}$ can always be attained. If true, this is best possible. In this section, we show that $V=[n]$ attains the equality and discuss further tight instances for the \LRP.

Let $V=[n]$. Recall Dirichlet's approximation theorem, whose proof follows easily from the \emph{Pigeonhole principle}: for every $t\in\RR$ and $n\in \NN$, there exist integers $i,j$ with $1\leq i\leq n$ such that 
$$
|ti-j|\leq \frac{1}{n+1}.
$$
Equivalently, for every $t\in \RR$, there is $i\in [n]$ with $\|ti\| \leq \frac{1}{n+1}$, and $\kappa(V)\leq \frac{1}{n+1}$. This proves the upper bound in~\eqref{eq:SODP}:
$$
\kappa(n)\leq \frac{1}{n+1}.
$$ 
By considering $t=\tfrac{1}{n+1}$, we obtain $\kappa(V)= \frac{1}{n+1}$. Figure~\ref{fig:plot2} illustrates the examples with two and three runners with these sets of speeds.

\begin{center}
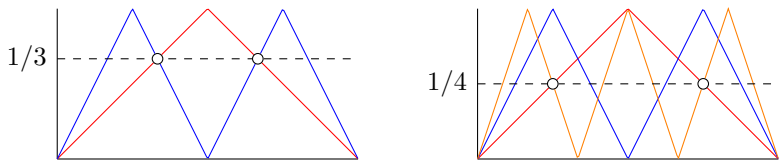
\begin{figure}[h]
\begin{center}
\begin{tikzpicture}[scale=4]
\draw (0,0)--(0,0.5) (0,0)--(1,0);
\draw [red, domain=0:1,samples=200] plot (\x, {min(\x-floor(\x),ceil(\x)-\x});
\draw [blue, domain=0:1,samples=200] plot (\x, {min((2*\x-floor(2*\x)),(ceil(2*\x)-2*\x)});
\draw[dashed] (0,1/3)--(1,1/3);
\foreach \i in {1/3,2/3}
{
\draw[fill,white] (\i,1/3) circle(0.5pt);
\draw (\i,1/3) circle(0.5pt);
}
\node[left] at (0,1/3) {\footnotesize{$1/3$}};
\end{tikzpicture}
\hspace{5mm}
\begin{tikzpicture}[scale=4]
\draw (0,0)--(0,0.5) (0,0)--(1,0);
\draw [red, domain=0:1,samples=200] plot (\x, {min(\x-floor(\x),ceil(\x)-\x});
\draw [blue, domain=0:1,samples=200] plot (\x, {min((2*\x-floor(2*\x)),(ceil(2*\x)-2*\x)});
\draw [orange, domain=0:1,samples=200] plot (\x, {min((3*\x-floor(3*\x)),(ceil(3*\x)-3*\x)});
\draw[dashed] (0,1/4)--(1,1/4);
\foreach \i in {1/4,3/4}
{
\draw[fill,white] (\i,1/4) circle(0.5pt);
\draw (\i,1/4) circle(0.5pt);
}
\node[left] at (0,1/4) {\footnotesize{$1/4$}};
\end{tikzpicture}
\end{center}
\caption{Plot of $\|tv_i\|$ for the sets of speeds $\{1,2\}$ (left) and $\{1,2,3\}$ (right).}\label{fig:plot2}
\end{figure}
\end{center}

{We say a set of speeds $V$, $|V|=n$, is a \emph{tight instance} for the \LRP\ if $\kappa(V)=\tfrac{1}{n+1}$. Since the \LRP\ is invariant by dilations, the sets $V=\{k,2k,\dots, nk\}$ are tight for any $k\in \NN$. These are the only tight instances for $n=1,2,3$} as shown by Cusick~\cite{C1973,C1974}. After a comment by P. Flor, Wills~\cite{W1968} identified the following tight instances:
\begin{align*}
V&=\{1,3,4,7\},\\
V&=\{1,3,4,5,9\},\\
V&=\{1,2,3,4,5,7,12\}.
\end{align*}
For $n=4$, Chen \cite{C1991} shows that there are no more instances, up to dilations (see also \cite{BS2008b}), and the same is true for $n=5$ as shown by Bohman, Holzman and Kleitman~\cite[Theorem 3]{BHK2001}.

Goddyn and Wong~\cite{GW2006} give the following tight sets:
\begin{align*} 
V&=\{1,4,5,6,7,11,13\},\\ 
V&=\{1,2,3,4,5,6,7,8,9,10,11,13,24\},\\
V&=\{1,2,3,4,5,6,7,8,9,10,11,12,13,14,15,16,17,19,36\},
\end{align*}
identifying a pattern: sometimes, accelerating a single speed $r\in [n]$ from $V=[n]$ produces a tight instance. This allows them to display infinite families of tight instances different from dilations of $[n]$. Precisely, for $n\ge 4$, replacing one integer $r\in [n]$ by a multiple $mr$ (accelerating runner $r$ by a factor $m$) leads to another tight instance if and only if $r$ has a common factor with all elements in the interval $[(n+1)-r, m(n+1-r)-1]$. For instance, if $r=n-1$ and $m=2$, we obtain that 
\begin{equation}\label{eq:SPEI}
    V=\{1,2,\dots, n-2, n,2(n-1)\},
\end{equation}
is a tight instance for every $n=6t+1$, $t\in \NN$. (Observe that the cases $t=1,2,3$ were already presented above.)

Moreover, accelerating simultaneously a set of runners from $[n]$ such that each individual one satisfies the above arithmetic property, leads again to a tight instance, leading to the following result.
\begin{theorem}[\cite{GW2006}]\label{thm:tightex} Let $m_1, m_2, \ldots ,m_n$ be positive integers. Suppose that each $r\in [n]$ has a common factor with every element in the interval $[(n+1-r), m_r(n+1-r)+1]$, then 
$$
V=\{m_1,2m_2,\ldots ,nm_n\}
$$
is a tight instance for the \LRC.
\end{theorem}

One can construct infinite families of sets $\{m_1, m_2,\ldots ,m_n\}$ for which the sufficient condition in the above theorem holds, thus providing additional infinite families of tight instances. 
This problem is closely related to the \emph{Jacobsthal problem}~\cite{J1961}, a notorious problem in number theory which asks for the length $g(r)$ of a largest interval of consecutive integers having a nontrivial common factor with $r$. Erd\H os~\cite{E1962} showed that there are infinitely many values of $r$ for which $g(r)\geq \log{r}$, which could be used to produce tight instances with $v_n=  2n - \Theta(\log{n})$. On the other hand, Pomerance~\cite{P2022} showed that there is a constant $c>0$ such that if $n<v_n<2n-c\log^2 n$, then $V$ is not tight; see Section~\ref{subsec:lin_bounded}. Bridging the gap would require new ideas on the Jacobsthal problem.

The problem of providing a complete characterization of tight instances is still widely open. In particular, the converse of Theorem~\ref{thm:tightex} does not hold in its full generality; see~\cite[Section~3]{GW2006}.

\section{Gap of loneliness}\label{sec:gap}

One approach to the \LRC\ is to derive lower bounds on the gap of loneliness $\kappa(n)$. As stated in~\eqref{eq:SODP}, Wills~\cite{W1967,W1968} observed that 
\begin{equation}\label{eq:EIRD}
\kappa(n) \geq \frac{1}{2n}.
\end{equation}
The proof follows from a simple argument, which we frame in the language of Probability Theory: Fix a set of speeds $V$ of size $n$. Consider the uniform probability space on $(0,1)$. For each $v\in V$, let $A_v$ be the event composed of the times $t\in (0,1)$ for which $\|tv\|< 1/2n$. Then $\Pr(A_v)< 1/n$. By the union bound, $\Pr(\cup_{v\in V} A_v)<1$ and there exists $t\in (0,1)$ where all the runners are at distance at least $1/2n$ from the origin.

Significant improvements on this trivial lower bound have remained elusive. Modest progress has been obtained in the literature, which we now discuss. Chen~\cite{C1994} showed that
$$
\kappa (n)\ge \frac{1}{2n-1-\frac{1}{2n-3}}.
$$
Moreover, Chen and Cusick~\cite{CC1999} proved the following partial improvement.
\begin{theorem}\label{thm:cc} Let $a$ and $n$ be positive integers. For every prime $p$ satisfying
$$
p\ge \max\{ 2a(n-1)-1, 2(a-1)n+2\}
$$
we have 
$$
\kappa (n) \ge \frac{a}{p}.
$$
\end{theorem}

It follows from Theorem~\ref{thm:cc} with the choices $a=1,2$ that if $2n-3$ is prime then $\kappa (n)\ge 1/(2n-3)$ and, if $4n-5$ is a prime, then $\kappa (n)\ge 2/(4n-5)$. 


One idea to improve the union bound obtained in~\eqref{eq:EIRD} is to evaluate more precisely the intersection between the sets in the family $\{A_v:\, v\in V\}$. By knowing only the measure of pairwise intersections, a Bonferroni-type inequality by Hunter~\cite{H1976} allowed the authors~\cite{PS2016} to get the improvement, for all $n$,
$$
\kappa (n)\ge \frac{1}{2n-2+o(1)}.
$$
This involves an exact calculation of $\Pr(A_{v_1}\cap A_{v_2})$, based on an analogous computation by Alon and Ruzsa~\cite{AR1999}. 

Tao~\cite{T2018} suggests an example for which the approach by approximating the Inclusion--Exclusion formula as above seems to be bound to a small improvement on the trivial bound, of the type $1/(2n)+O(1/n^2)$ as in the above results. {To this end, consider a set $V$ that contains all the primes $p_1,p_2,\dots, p_s$ between $n/4$ and $n/2$; by the Prime Number Theorem, $s\sim \frac{n}{4\log n}$. By exploiting the coprimality, one can show that
$$
\Pr(\cup_{i=1}^s A_{p_i}) = \big(1+O\big(\tfrac{1}{n}\big)\big) \sum_{i=1}^s \Pr(A_{p_i}).
$$
While this example shows that there is a set of size $s$ where the union bound has an error factor of order $O(1/n)$, it does not take into account the interaction produced by the remaining $n-s$ sets, for which intersections are likely to be unavoidable.}

{Tao approaches the problem from the perspective of \emph{Bohr Sets} and \emph{Generalised Arithmetic Progressions}, which are central topics in Additive Combinatorics. By obtaining a cruder estimation on $\Pr(A_{v_1}\cap A_{v_2}\cap A_{v_3})$, which correspond to Bohr sets of dimension $3$,} the following almost logarithmic improvement is obtained in~\cite{T2018},
\begin{equation}\label{eq:tao}
    \kappa (n)\ge \frac{1}{2n}+\frac{c\log n}{n^2 \log\log n}.
\end{equation}
The approach in~\cite{PS2016} also gives, 
$$
\kappa (V)\ge \frac{1}{2n}+ \frac{c}{n^2}\Big(\sum_{i=2}^n\frac{1}{v_i}\Big),
$$
which provides the logarithmic improvement if the sum of the velocity reciprocals is logarithmic. This includes the case where all speeds are linearly bounded; see Subsection~\ref{subsec:lin_bounded}.

\section{Small number of runners}\label{sec:small}

As mentioned in Section~\ref{sec:connections}, different approaches to the \LRP\ have provided in particular the solution of the \LRC\ for small values of $n$. Currently the conjecture is proved for all $n\le 6$ (namely, up to seven runners in Conjecture~\ref{conj:LRC_original}, the original formulation). We next review the different approaches used to solve the conjecture for these small values of $n$.

\subsection{Case $n=2$}
We need to prove that 
\begin{equation}\label{eq:n=2}
\kappa(2)=1/3.
\end{equation}
Several arguments are available. Let $V=\{v_1,v_2\}$, where $v_1<v_2$ are positive integers with $(v_1,v_2)=1$. A first proof of~\eqref{eq:n=2} was given by Wills~\cite{W1968} based on Bezout's identity. 

From the view-obstruction approach, as illustrated in Figure~\ref{fig:view} with $x=v_1$ and $y=v_2$, the line $y=2x$ is tangent to the corners of the translated squares with dilation $1/3$ and, for any smaller dilation of the inner square these lines miss all the squares. All other rays clearly hit the squares, showing~\eqref{eq:n=2}. (See e.g. Cusick~\cite[Introduction]{C1973}.)

Another simple way to deduce~\eqref{eq:n=2}, presented in~\cite{BGGST1998}, uses the observation made in Section~\ref{sec:reductions} which states that the time when the runners are furthest from the origin is of the form $t_0=\tfrac{\ell}{v_1+v_2}$ for some $\ell\in \NN$ with $\ell\leq \lfloor\frac{v_1+v_2}{2}\rfloor$. By solving the congruence $v_1t_0= \lfloor\frac{v_1+v_2}{2}\rfloor \mod (v_1+v_2)$, we obtain that 
$$
\kappa (V) = \frac{\lfloor(v_1+v_2)/2\rfloor}{(v_1+v_2)}\geq \frac{1}{3}.
$$
The equality is attained if and only if $V=\{k,2k\}$ for $k\in \NN$, showing that these are the only tight instances for $n=2$. Incidentally, for this simple case, the value of $\kappa (V)$ equals the maximum density of a set avoiding the two differences in $V$, the Motzkin problem mentioned in Section~\ref{subsec:chrom}, see e.g.~\cite[Example 1]{L2008}.

A fourth approach is given in~\cite{CP1984} using the discrete version of the \LRP\ discussed at the end of Section~\ref{subsec:chrom}. If both $v_1$ and $v_2$ are coprime with $3$ we readily have $\kappa_3(V)\ge 1$ and $\kappa (V)\ge 1/3$. Otherwise, if $m$ is the largest power of $3$ dividing $v_2$, say, then choose $N=3^{m+1}$ so that $\|v_2\|_N\ge N/3$. Then $\lambda_k=(1+3^mk)$ for $k\in \{0,\pm1,\pm 2\}$ satisfies $\lambda_kv_2=v_2 {\mod N}$  and at least one of the choices for $k$ satisfies $\|\lambda_kv_1\|_N\ge N/3$, leading again to~\eqref{eq:n=2}.

\subsection{Case $n=3$}

In this case, the conjecture states that
\begin{equation}\label{eq:n=3}
\kappa(3)=1/4.
\end{equation}
Betke and Wills~\cite{BW1972} give a proof from the perspective of Diophantine approximation, although they mention that their method does not extend to $n\ge 4$. 

Shortly after, Cusick provided two proofs of this case. In~\cite{C1973}, he gives a proof based on the view-obstruction paradigm: project all the inner cubes in $\RR^3_{>0}$ onto the faces of the unit cube supported on the planes $x=1$, $y=1$ and $z=1$, along the line of sight. By showing that this projection fully covers these three faces, one concludes~\eqref{eq:n=3}.

In~\cite{C1974} Cusick reproves the cases $n=2,3$ using integral formulae. For $v\in \NN$, let $F_{v}(x)$ be the indicator function that $\|v x\|> 1/4$. He shows
$$
I(V)=\int_{0}^1 F_{v_1}(x)F_{v_2}(x)F_{v_3}(x)\,dx = 0,
$$ only if $V=\{v_1,v_2,v_3\}=\{k,2k,3k\}$, for which the conjecture holds (with equality). Observe that if $I(V)>0$, the set of times where the origin is isolated has positive measure, thus proving~\eqref{eq:n=3}. Additionally, the method of proof yields the uniqueness of tight instances for $n=3$.

Cusick~\cite{C1982} later gives an elementary proof of~\eqref{eq:n=3} which is close in spirit to the discrete version argument given above for the case $n=2$, and he suggests that the argument could be extended to $n\geq 4$.

\subsection{Case $n=4$}
Elaborating on the method introduced in~\cite{C1982}, Cusick and Pomerance~\cite{CP1984} prove that
\begin{equation}\label{eq:n=4}
\kappa(4)=1/5.
\end{equation}
The proof involves computer aid to check some small cases and it is unclear whether it can be extended for $n\geq 5$.

Chen~\cite{C1991} gives an elementary proof based on case analysis for  $n\le 4$. Indeed, he proves another statement that is proved to be equivalent to the \LRC: for every $n$ and $n$-set $V$ of positive real numbers, there exists integers $t_1,t_2,\dots, t_n$ such that
$$
t_jv_i-t_iv_j\leq \frac{n}{n+1}v_j-\frac{1}{n+1}v_i \quad \text{for }i,j\in [n].
$$

A third simpler proof for $n=4$ is given by  Bienia, Goddyn, Gvozdjak, Seb\H o, and Tarsi~\cite{BGGST1998}, following the same ideas displayed in the third proof for the case $n=2$. They illustrate their method by giving yet another proof for $n=3$. In their argument for $n=4$, it is crucial that $n+1$ is a prime number.

\subsection{Case $n=5$}
The case $n=5$ is first proved by Bohman, Holzman and Kleitman~\cite{BHK2001}. Roughly speaking they separate the two faster runners from the three slower ones and show that the intervals of time where the two faster ones are not lonely can not cover the intervals where the three slower ones are lonely. The proof is framed in a nice geometric setting: hitting a polygon in a torus by a discrete group. As in most of the proofs, a case analysis of critical cases seems to prevent the method to be extended to larger values of $n$. The method identifies the two families of tight instances for this case: $\{k,2k,3k,4k,5k\}$ and $\{k,3k,4k,5k,9k\}$ for $k\in \NN$. 

A simpler proof for $n=5$ is given by Renault~\cite{R2004} which takes into account the distinct congruence classes of the speeds modulo six, studying them by case analysis. He also provides proofs with the same ideas for smaller values of $n$, including a compact proof for $n=4$.

\subsection{Case $n=6$}
The currently largest value of $n$ for which the conjecture is proved is $n=6$, by Barajas and Serra~\cite{BS2008}. Here again the proof takes advantage of the fact that $n+1$ is prime and the analysis is made through the distinct congruence classes of the speeds modulo seven. A lemma which was implicitly used for the case $n=4$ in~\cite{CP1984} is explicitly formulated (and called the Prime Filtering Lemma). This Lemma  easily allows one to move to the most intricate case where five out of the six speeds are coprime with $7$. The technique could be used for larger $n$ such that $n+1$ is prime, but already for the case $n=10$ the amount of case analysis becomes cumbersome.

\section{Special families}\label{sec:special}

The conjecture has been proved for several special families of sets of speeds. In some cases, even the exact value of $\kappa (V)$ is known.

\subsection{Simple remarks and further reductions} \label{sec:simple_rem}

In Section~\ref{sec:reductions} we have discussed some reductions to the original conjecture. There are some additional simple observations which make the problem easier or even trivial under additional hypotheses. 

A first observation is that if no number in $V$ is a multiple of some $k\in [n+1]$ then
$$
\kappa (V)\ge \frac{1}{k}\kappa_{k} (V)= \frac{1}{k}\ge \frac{1}{n+1},
$$
so that $V$ satisfies the conjecture.

Another simple observation is that if there is a natural $N$ such that $V\subset [N/{(n+1)}, N-N/(n+1)]$ then we plainly have 
$$
\kappa (V)\ge \frac{1}{N}\kappa_N(V)\ge \frac{1}{n+1}.
$$
The above condition can be simply written as $v_n\le nv_1$ (recall that speeds are assumed to be increasing), in which case $N=v_1+v_n$ is a suitable value. This means that, for every set $V$, every  translate $a+V$ with $a\ge (v_n-nv_1)/(n-1)$ satisfies 
\begin{equation}\label{eq:VNSI}
\kappa (a+V)\ge \frac{1}{n+1}.    
\end{equation}

By using the polytopal perspective discussed in Section~\ref{subsec:zonotope}, Beck, Hosten and Schymura~\cite{BHS2019} proved  that the conditions $v_{n-1}\le (n-2)v_1$, or $v_{n-2}\le (n-2)v_1$ and $v_{n-1}\ge nv_{n-2}$, also ensure $\kappa (V)\ge 1/(n+1)$. Additional conditions of this flavour are also obtained by Bhardwaj, Narayanan and Venkataraman~\cite{BNV2022} with a similar approach.

Finally, adding a sufficiently large speed $v_n$ to a set $V'$ with $n-1$ speeds for which the conjecture is satisfied, i.e. $\kappa (V')\ge \tfrac{1}{n}$, results in a set $V=V'\cup \{v_n\}$ which satisfies the $\LRC$. Suppose that $t_0$ is a time where $V'$ is at distance at least $1/n$ from the origin. Then, there exists a time interval centered at $t_0$ and of length
$$
\ell=\frac{2}{v_{n-1}}\left(\frac{1}{n}-\frac{1}{n+1}\right)=\frac{2}{n(n+1) v_{n-1}},
$$
such that all runners are still at distance at least $1/(n+1)$ from the origin. If $v_n\geq  nv_{n-1}$, during this time the fastest runner covers a piece of size 
$$
v_n \ell \geq  \frac{2}{n+1},
$$
and thus there is a time in the interval where its distance from the origin is at least $1/(n+1)$.

\subsection{Sets related to arithmetic progressions}\label{sec:param}

A natural case to analyse is a modification of the tight instance $V=[n]$. Liu~\cite{L2008} obtained the exact value for $\kappa(V)$, where $V$ is obtained from $[n]$ by removing an arithmetic progression.

\begin{theorem} Let $n,k,s$ be positive integers with $n\ge (s+1)k$. Then, for 
$$
V=[n]\setminus \{k,2k,\ldots ,sk\},
$$
we have 
$$
\kappa (V)=\begin{cases} 
\frac{1}{2k}, & n=(s+2)k-1 \; \text{and}\; s \; \text{even};  \\
\frac{s+1}{m+sk+t}, & \text{otherwise}.
\end{cases}
$$
where $t$ is the smallest positive integer such that $\text{gcd}(n+t,k)=1$.
\end{theorem}

In particular, removing one element $k\le n/2$ from $[n]$ results in a set $V$ for which $\kappa (V)>1/|V|$  (except for $V=\{1,3,4\}$ for which $\kappa (V)=2/7$), a stronger conclusion than the conjectured one. The arithmetic progression structure of the removed set from $[n]$ may be relevant, as some of the sporadic tight instances show. 

Another variation is to remove a full interval from $[n]$. Let $1\le a<b\le n$ and set $V=[n]\setminus [a,b]$. The exact value of $\kappa (V)$ is obtained by Liu~\cite{L2008} building upon previous results on the chromatic numbers of distance graphs by Wu and Lin~\cite{WL2004}.

\begin{theorem} Let $2\le a,b\le n$ be positive integers with $a+1\le b\le 2a-1$ and set $V=[n]\setminus [a,b]$. Then
$$
\kappa (V)=\begin{cases}
\frac{2}{n+1}, & 2a\le n<2b;\\
\frac{2}{n+a+t}, & n\ge 2b,
\end{cases}
$$
where $t$ is the smallest positive integer such that $\gcd (n+a+t,y)=1$ for some $a\le y\le \min\{a+t-1,b\}$.
\end{theorem}

Pandey~\cite{P2010} generalises $V=[n]$ by considering the case of arithmetic progressions and obtains the following lower bound.
\begin{theorem} Let $V=\{a,a+d,\cdots ,a+(n-1)d\}$ for some positive integers $a,d$ and $n\ge 2$. Then
$$
\kappa (V)\ge 
\begin{cases}
\frac{1}{2}, & d \; \text{even;}\\
\frac{1}{2}\frac{2a+(n-1)(d-1)}{2a+(n-1)d},& d \; \text{odd.}
\end{cases}
$$
\end{theorem}

The bound is good enough not only to prove the Conjecture for arithmetic progressions but also for dense subsets of arithmetic progressions.

\begin{corollary}\label{cor:ap} Let $V$ be an $n$-subset of an arithmetic progression $P=\{a,a+d,\ldots ,a+(k-1)d\}$ of length  $k\le 2n-3$. If $a$ and $d$ are not both one, then
$\kappa (V)\ge \tfrac{1}{n+1}$.
\end{corollary}

\subsection{Lacunary sequences}\label{sec:lacunary}


A sequence $(v_n)_{n\ge 1}$ is \emph{$\epsilon$-lacunary} (or simply \emph{lacunary}) if there is $\epsilon>0$ such that $v_{n+1}\ge (1+\epsilon)v_{n}$ for all $n\ge 1$. 

Rusza, Tuza and Voigt~\cite{RTV2002} and Pandey~\cite{P2009} studied lacunary sequences where each speed is roughly at least twice the previous one. The following result was obtained by Barajas and Serra.
\begin{theorem}[\cite{BS2009}]\label{thm:pandey_lacu} For every $n\geq 2$ and every $n$-set $V$ satisfying
$$
v_{i+1}\ge 2 v_i \quad \text{for }i\in [n-1],
$$
we have $\kappa(V)\geq \frac{1}{n+1}$.
 \end{theorem}
Therefore the \LRC\ holds for a quantitatively specified version of lacunary sequences, in particular, this result covers sequences that grow at least exponentially.


A classical result by Weyl~\cite{W1916} implies that  {if $(v_n)_{n\ge 1}$ is lacunary then} the sequence $({tv_n \mod 1})_{n\geq 1}$ is uniformly distributed in $[0,1)$ for almost all $t\in \RR$. Motivated by problems on the chromatic number of distance graphs, Erd\H{o}s~\cite{E1975} asked if  {for any $\epsilon>0$ there is $\delta>0$ such that for any $\epsilon$-lacunary sequence $(v_n)_{n\ge 1}$ 
 \begin{equation}\label{eq:SPAS}
 \sup_{t\in (0,1)}\inf_{n\ge 1} \|tv_n\|\ge \delta. 
 \end{equation}
This was proved independently by Mathan~\cite{M1980} and Pollington~\cite{P1979}. Since then there have been a number of results improving the dependence of $\delta$ on $\epsilon$~\cite{K2001,RTV2002}, leading to the following result of Peres and Schlag~\cite{PS2010}}: 
  $$
  \sup_{t\in (0,1)}\inf_{n\ge 1} \|tv_n\|\ge c\epsilon |\log \epsilon|^{-1},
  $$ 
  for some absolute constant $c>0$.  If $E_i=\{t\in (0,1):\|tv_i\|<\frac{1}{240 M\log M}\}$ for $M=\lceil 1/\epsilon\rceil$ then it follows from a directed lopsided version of the Local Lov\'asz Lemma that $\big(\!\cap_{i\ge 1} \overline{E_i}\big)\neq \emptyset$, where $\overline{E_i}$ is the complement of $E_i$, proving the result (actually in more generality which includes finite unions of lacunary sequences).  
  
The parameter in the RHS of~\eqref{eq:SPAS} can be understood as an infinite version of $\kappa(V)$ as defined in~\eqref{eq:xi_kappaV}, which directly relates it with the \LRP. Dubickas~\cite{D2011} used the same approach to improve Theorem~\ref{thm:pandey_lacu} for large $n$.
  
 \begin{theorem}\label{thm:dubickas} There is a positive integer $n_0$ such that, for every $n\ge n_0$ and every $n$-set $V$ satisfying 
 $$
 v_{i+1}\ge \left(1+\frac{22\log n}{n}\right)v_i \quad \text{for }i\in [n-1],
 $$
 we have $\kappa(V)\geq \frac{1}{n+1}$.
 \end{theorem}
 
In fact, Dubickas~\cite{D2011} proves a slightly stronger result that only needs the ratio of speeds that are far apart to be bounded away. Precisely, it suffices that $n\geq 32$ and
\begin{equation}\label{eq:SPDS}
v_{i+\lceil\frac{n+1}{12e}\rceil}\geq (n+1) v_i  \quad \text{for }i=1,2,\dots, n-\lceil\tfrac{n+1}{12e}\rceil.    
\end{equation}

In contrast with Theorem~\ref{thm:pandey_lacu}, the sequences for which this theorem holds have polynomial growth, of the order of $n^{22}$. With the constant $22$ replaced by $33$, the theorem holds for $n\ge 16342$, which gives an idea of the size of $n_0$ in Theorem \ref{thm:dubickas}. 

The condition in~\eqref{eq:SPDS} can be seen as a sort of counterpart of the results presented in Section~\ref{sec:simple_rem}, where the \LRC\ is proved in the cases where $v_n$ is not too large with respect to $v_1$, e.g. $v_n\leq nv_1$.

A refinement of Theorem~\ref{thm:dubickas} using the Local Cut Lemma of Bernshteyn~\cite{B2017} was obtained by Czerwi\'{n}ski~\cite{C2018} by weakening the lacunary property to hold only for the fastest runners and allowing the slower ones to take any value.

\section{Bounded speeds}\label{sec:bounded}

In this section we discuss the \LRP\ for families of sets of speeds that are bounded in some sense.

\subsection{Reduction to bounded speeds}\label{sec:bounded1}

A remarkable result by Tao~\cite{T2018} is that it suffices to check that $\kappa (V)\ge \frac{1}{n+1}$ for a finite number of $n$-sets $V$ to conclude that the conjecture holds 	 {for a given $n$}. The statement can be formulated as follows.


{
\begin{theorem}\label{thm:bounded} There is a computable constant $C>0$ such that for every $n_0$, if $\kappa (V)\ge \frac{1}{n+1}$ for all $n$-sets $V$ with $n\le n_0$ and $\max (V)\le n^{Cn^2}$, then
$$
\kappa (n)\ge \frac{1}{n+1} \quad\text{for all }n\leq n_0.
$$
\end{theorem}
}

The proof is based on an appropriate embedding of $V$ in a  proper centered $r$--dimensional arithmetic progression $Q$ of size $|Q|=O(n^{Cn^2})$ with $r\leq n$. Such a generalized arithmetic progression $Q\supset V$ is of the form 
$$
Q=\{n_1w_1+\cdots +n_rw_r: n_i\in \ZZ, |n_i|\le N_i \text{ for all }i\in [r]\}
$$
for some integers $w_1,\ldots ,w_r$ and $N_1,\ldots ,N_r$. If $r=1$ then in this embedding all speeds lie in a centered arithmetic progression and, by dividing them by difference $w_1$, one obtains the result. It is worth noticing that this case captures all the tight instances of the \LRC. If $r\ge 2$ then, one can project $V$ into a set $V'$ in such a way that the projections of two distinct speeds coincide, so $|V'|\leq n-1$. By induction, for the set of speeds $V'$, there exists a time where every runner is at distance at least $1/n$ from the origin. A delicate argument allows then to transfer the result from $V'$ to $V$ only losing a small additive factor $1/5n^2$ on the gap, which is enough to prove the lower bound $1/(n+1)$. In this argument it is key to have good upper bounds on the size of Bohr sets, i.e. the sets of times when all runners are close to the origin.

The statement in Theorem \ref{thm:bounded} has been recently reproved with the quantitative bound significantly improved by Malikiosis, Santos and Schymura \cite{MSS2025}.

\begin{theorem}[\cite{MSS2025}]\label{thm:bounded+} Suppose that the LRC holds for $n$ runners. Then the conjecture holds for every set  of $n+1$ runners with velocities $V=\{0,v_1,\ldots ,v_n\}$ satisfying $\gcd(v_1,\ldots ,v_n)=1$ and 
$$
v_1+\cdots +v_n>\binom{n+1}{2}^{n-1}.
$$
\end{theorem}

The proof of Theorem \ref{thm:bounded+} uses the geometric approach of the LRC in terms  of zonotopes which has been   discussed in Subsection \ref{subsec:zonotope}. With this geometric interpretation all what is needed for the proof of Theorem \ref{thm:bounded+} is the classical Minkowski  theorem on the first successive minima of a convex body, which allows one to reduce the problem to one dimension less as long as the number of lattice points of the zonotope is at least $\binom{n+1}{2}^{n-1}$, this being  the real source of the lower bound in the theorem. One additional advantage of this zonotope  approach is that it provides an analogous result for the so--called shifted version of the Lonely Runner Problem, see Subsection \ref{subsec:shifted}.

\subsection{Linearly bounded speeds}\label{subsec:lin_bounded}

Even if the bound  on the interval where the sets $V$ can be taken is not practical from the computational point of view, Theorem~\ref{thm:bounded} prompts the question of exploring the \LRC\ for sets contained in bounded intervals. In this direction the following is also proved in~\cite{T2018}.

\begin{theorem}\label{thm:bounded2} For every $n\geq 2$ and every $n$-set $V$, if $v_n<1.2 n$ then $\kappa (V)\ge \frac{1}{n+1}$.
\end{theorem} 

Tao poses the question of whether the constant 1.2 can be improved to $2$, which would include tight instances different from multiples of $[n]$; see Section~\ref{sec:tight}. Prompted by this question, Bohman and Peng~\cite{BP2022} develop the approach of coprime mappings presented in Section~\ref{subsec:coprime_map}. As an application of Theorem~\ref{thm:cop}, they obtained the following approximate solution to the question.

\begin{theorem}\label{thm:bounded3}  There exists a constant $c>0$ such that for sufficiently large $n$ and every $n$-set $V$, if $n<v_n\le 2n-\exp(c(\log\log n)^2)$, then $\kappa(V)> \tfrac{1}{n+1}$. 
\end{theorem}
Shortly after, Pomerance~\cite{P2022} weakened the required lower bound on $m$ in Theorem~\ref{thm:cop}, which implies that the error term in the upper bound of $v_n$ in Theorem~\ref{thm:bounded3} can be replaced by $c(\log n)^2$. Further improvements on this would require improving upon the best known bounds for the Jacobsthal problem; see Section~\ref{sec:tight}. 

One should compare these results with Corollary~\ref{cor:ap} by Pandey~\cite{P2010}. A simple consequence of it is that the $\LRC$ holds provided that $v_n\le 2n-3$ except in the case where $v_1=1$ and $V$ is not included in an arithmetic progression with difference at least two.

\section{Random runners}\label{sec:random}

{A natural question is to ask if the \LRC\ holds for random subsets of speeds. If we consider random real speeds $v_1,v_2,\dots, v_n$, for any reasonable notion of randomness, they will be almost surely linearly independent over $Q$. By Kronecker's theorem (Theorem~\ref{thm:kron}) with $a_i=1/2$ for all $i\in [n]$,
$$
\| tv_i\| \geq \frac{1}{2}-\epsilon, \quad \text{for every $i\in [n]$}\;.
$$
}

As the problem can be reduced to positive integer speeds, it is natural to ask if the previous conclusion still holds true for random integer speeds. For a fixed $n$, Czerwi\'{n}ski~\cite{C2012} showed that asymptotically almost surely as $N$ tends to infinity, an $n$-subset $V$ uniformly chosen from the interval $[N]$ also satisfies the conjecture in a very strong sense.

\begin{theorem}\label{thm:random} Let $n\ge 2$ be an integer. For every $\epsilon>0$, the probability that a subset  $V$ of size $n$ chosen uniformly at random among the $n$--subsets of $[N]$ satisfies
$$
\kappa (V)\ge \frac{1}{2}-\epsilon,
$$
tends to one as $N\to \infty$.
\end{theorem}
The proof shows that almost all subsets of $\ZZ_p$, $p$ a prime, satisfy a property called $L$--independence (roughly speaking, no solutions to homogeneous linear equations with coefficients bounded by $L$), and for such sets, Fourier analytic techniques show that they satisfy the inequality stated in the above Theorem. Kravitz~\cite{K2021} noticed that the same argument can be applied to random runners extracted from an arbitrary set $S\subset \NN$ of cardinality $N$. 

{The previous discussion hints that, as already observed in Section~\ref{sec:tight}, the instances for which the \LRP\ is hard are well-structured.}

From the perspective of the chromatic number of  distance graphs  discussed in Section~\ref{subsec:chrom}, this in particular shows that asymptotically almost surely the chromatic number of distance graphs is at most three. More precisely, the argument in~\cite{C2012} shows that the chromatic number of the circulant graph $Cay (\ZZ/p\ZZ,V)$ for $|V|\le (\log p/\log\log p)^{1/2}$ is $3$ with probability tending to one as $p$ tends to infinity. Alon~\cite{A2013} addressed the more general question of studying the chromatic number of random Cayley graphs. In particular he obtains, by a probabilistic argument, a quantitative improvement on the bounds implicit in the statement of Theorem~\ref{thm:random}. One way to state the result relevant to the present discussion~\cite[Theorem 3.7]{A2013} is as follows.

\begin{theorem} Let $N$ be a positive integer and $\epsilon>0$. For $n\le (1-\epsilon)\log_3 N$, almost all sets $V\subset [N]$ with $|V|=n$ satisfy
$$
\kappa (V)\ge \frac{1}{3}.
$$
\end{theorem}

\section{Variations and generalizations}\label{sec:variations}

In this section we discuss several variations and generalizations on the \LRP\ which have appeared in the literature.

\subsection{Invisible runners}

As it has been mentioned in Section~\ref{sec:special}, the removal of one element in $[n]$ results in a set which satisfies the Conjecture. Czerwi\'{n}ski and Grytczuk~\cite{CG2008} proved that this is a general phenomenon which they call the {\it invisible runner} problem.

\begin{theorem}\label{thm:CG}
For every $n\in \NN$, every $n$-set $V$ of positive integers contains an element $v\in V$ such that
$$
\kappa (V\setminus \{v\})\ge \frac{1}{n}.
$$
\end{theorem}

We provide a simple proof of this result: Assume that $\kappa(V)< \tfrac{1}{n}$, then the average number of runners at distance $\kappa$ from the origin is $2\kappa n< 2$. So there exists a time where only one runner, say $v\in V$, is at distance less than $\tfrac{1}{n}$ from the origin. This time certifies that $\kappa (V\setminus \{v\})\ge \frac{1}{n}$.

Actually, if $s\ge 1$ runners are allowed to be invisible, then the lower bound for $\kappa$ can be improved: every set $V$ of $n$ positive integers  contains a subset $S\subset V$ with $|S|=s$ such that $\kappa (V\setminus S)\ge (s+1)/2n$.   

{In~\cite{PS2016} the authors slightly extend Theorem~\ref{thm:CG} above by showing that, for a set $V$, either the $\LRC$ holds or there are at least four speeds $v\in V$ such that $\kappa (V\setminus \{v\})\geq 1/n$. This result is proved using the notion of dynamic interval graphs. 
}

{
In relation to the problem of invisible runners, a weakening of the original conjecture was asked by Joel Spencer\footnote{Communicated to the authors by Jarek Grytczuk.}.
\begin{conjecture}[Single \LRC]\label{conj:LRC_single}
For every $n\in \NN$ and every $n$-set $V$ of pairwise distinct real numbers, there exists $i\in [n]$ and $t\in
\RR$ such that
\begin{equation}\label{eq:ALSD2}
\min_{j\neq i} \|t(v_j -v_i)\|\geq \frac{1}{n}.
\end{equation}
\end{conjecture}
}

\subsection{Spectrum of the loneliness gap}

{Kravitz~\cite{K2021} focused on the study of \emph{tight} and \emph{almost tight} instances. Based on his results, he posed the conjecture that, if a set of speeds is not almost tight, then its gap of loneliness will be uniformly bounded away from $\tfrac{1}{n+1}$.} More precisely, they conjectured that for every $n\in \NN$ and every $n$-set $V$ of positive integers we either have
$$
\kappa (V)= \frac{s}{sn+1} \; \text{for some } s\in \NN, \quad \text{(almost tight instances)} 
$$
or $\kappa (V)\ge \frac{1}{n}$.

Thus, according to Kravitz's Conjecture, the spectrum of possible values of $\kappa (V)$ in the interval $[\frac{1}{n+1},\frac{1}{n})$ is contained in the finite set 
$$
T=\left\{\frac{1}{n+1},\frac{1}{n+1/2},\ldots,  \frac{1}{n+1/n}\right\}.
$$

The conjecture is motivated by the fact that any value $\tfrac{s}{sn+1}$ in $T$ is reached by the set $\{1,2,\ldots ,n-1, ns\}$ and it asserts that no other value is possible. Kravitz proved the conjecture for $n=2,3$ and provided partial support for $n=4,6$ proving it whenever the fastest runner is much faster than the second fastest one: $v_n\ge 4v_{n-1}^4$. In the same paper, Conjecture~\ref{conj:spec} is also proved for sets with bounded speeds, specifically whenever $v_n\le 1.5n$ (see also Section~\ref{subsec:lin_bounded}). 

Fan and Sun~\cite{FS2023} disproved Kravitz's conjecture for $n=4$ by showing that $\kappa(3,8, 11, 19) = 7/30\in (3/13,4/17)$; a similar counterexample is given for $n=6$. They propose the following weakening of   the conjecture.

\begin{conjecture}[(Amended) Loneliness Spectrum Conjecture]\label{conj:spec} For every $n\in \NN$ and every $n$-set $V$ of positive integers we either have
$$
\kappa (V)= \frac{s}{sn+k} \; \text{for some } s,k\in \NN,
$$
or $\kappa (V)\ge \frac{1}{n}$.
\end{conjecture}

This line of research has been further explored by Giri and Kravitz \cite{GK2025}, by studying the spectra in higher dimensional tori. As an illustration of the applications of the approach they provide an improvement of Theorem \ref{thm:bounded} on the bounded speeds which gives a weaker bound than the one obtained in \cite{MSS2025}.  Jain and Kravitz~\cite{JK2024} have studied the Lonely Runner spectrum relative to $2$-dimensional
subtori $U\subset (\RR\backslash \ZZ)^n$, showing that it has a very rigid arithmetic structure. They used it to characterize the spectrum in different intervals for $n=3,4,6$.

\subsection{Shifted \LRC}\label{subsec:shifted}

The \LRC\ assumes that all runners start at the origin. From the perspective of the \LRC\ stemming from billiard ball trajectories (see Section~\ref{subsec:VO_Bill}), a natural generalization is to relax this assumption and let each runner start at a different position on the circular track; this is known as the \emph{Shifted \LRP}. The corresponding conjecture, which postulates the same bound in this case is explicitly formulated by Beck, Hosten and Schymura~\cite{BHS2019}, and its origin is attributed to Wills. 

In the shifted version it is important to insist that the speeds are pairwise distinct, in contrast to Conjecture~\ref{conj:LRC}. If the condition on distinct speeds is lifted, then the trivial lower bound $\tfrac{1}{2n}$ (see Section~\ref{sec:tight}) becomes tight, as witnessed by the set of speeds all equal to one and starting points $\{0,1/n,\cdots, (n-1)/n\}$. The lower bound in this case was found by Schoenberg~\cite{S1976} in the equivalent context of billiard trajectories; see also~\cite{BHS2019}. 

The \emph{covering radius} $\mu_K$ of a convex body $K\subset \RR^m$ is the minimum value of $\mu>0$ for which $\mu K+\ZZ^m$ covers $\RR^m$. Equivalently, $\mu_K$ is the minimum value of $\mu$ such that every translate of $\mu K$ meets a point in $\ZZ^m$, see e.g. \cite{GK1987}. Following the zonotopal interpretation of the \LRP\  described in Subsection \ref{subsec:zonotope}, the shifted  \LRP\ can be phrased in terms of covering radii as following: for every lattice zonotope $Z$ generated by $n$ vectors in general position in $\ZZ^{n-1}$, 
$$
\mu_{Z-\textbf{x}}=\tfrac{n-1}{n+1}
$$.

The case $n=2$ of the shifted version of the conjecture was proved in~\cite{BHS2019}. The case $n=3$ has also been confirmed by  Cslovjecsek, Malikiosis, Nasz\'odi and Schymura~\cite{CMMS2022}. Both proofs use the perspective of the covering radius of zonotopes. An additional proof for $n= 3$ has been obtained independently by Rifford~\cite{R2022}.

Yet another proof of the case $n=3$ can be found in Malikiosis, Santos and Schymura \cite{MSS2025}, where the covering radius of zonotopes is also used to show that the shifted conjecture will follow for $n=4$ if it can be proved for all sets of speeds $V$ with $v_n\le 195$. This is a refinement of a general result which states that Theorem \ref{thm:bounded+} holds for the shifted version of the \LRP\ conditional to an open problem in 2-dimensional geometry named as the \emph{Lonely Vector Problem}. An explicit computation is performed by Alcantara, Criado and Santos \cite{ACS2025} verifying the case $n=4$  of the shifted conjecture.

\subsection{Time to get lonely}

Given that the \LRC\ holds, it is natural to ask what is the smallest $t\in (0,1)$ for which the origin is lonely, which we denote by $t_0$\footnote{This question had been already considered by Chen in the context of a set of speeds that are pairwise linearly independent over $\QQ$~\cite{C2000}.}. 
The reduction on times given in Section~\ref{sec:reductions} does not provide any meaningful upper bound, as the bound depends on the speed set. As the \LRP\ is invariant by dilations, we renormalise the time by the slowest speed
$$
\hat{t}_0 = t_0v_1.
$$ 
Rifford~\cite{R2022} formulated the following strengthening of the \LRC.
%
%

{\begin{conjecture}[Timely \LRC]\label{con:LRC_time} 
For every $n\in \NN$ there is $N$ such that for every $n$-set $V$ of positive speeds
$$
\hat{t}_0\le N\;.
$$
\end{conjecture}
In words, it says that there is a uniform bound $N$ only depending on the dimension, such that the origin will be made lonely before the slowest runner has completed $N$ laps. In particular, the conjecture aligns with Theorem~\ref{thm:bounded}, which states that the \LRP\ can be reduced to a finite number of instances. We believe that a detailed analysis of the proof of the latter could prove Conjecture~\ref{con:LRC_time}, conditional on the existence of such $\hat{t}_0$.}

In~\cite{R2022} the conjecture is proved for $n\le 5$.

Bhardwaj, Narayanan and Venkataraman~\cite{BNV2022}, took a different approach and studied at which times the origin is lonely, posing the following conjecture, which is supported by some simulations.
\begin{conjecture} For every $n\in \NN$ and every $n$-set $V$ of positive integer speeds there is a positive integer $M$ such that if
$$
t= \frac{M}{2^{\lceil \ln_2 v_n +1\rceil}(n+1) v_n},
$$ 
then $\min_{v\in V} \| tv\|\geq \tfrac{1}{n+1}$.
\end{conjecture}

\subsection{Lonely Rabbit Problem}

Motivated by problems in Diophantine approximation, a variation on the \LRP\ modifies the definition of $\kappa (V)$ as follows. In Section~\ref{subsec:diophantine} we defined  
\begin{equation*}
\begin{aligned}
    \kappa(V)&=\sup_{t\in (0,1)}\min_{v\in V}\|tv\|,\\
    \xi(V)&=\sup_{t\in \ZZ}\min_{v\in V}\|tv\|.
\end{aligned}    
\end{equation*}
While in the first case we have runners that move continuously through the unit interval, in the second one, the runners become rabbits jumping simultaneously at discrete times with hop lengths given by $V$.  In this case, if $v$ is integer, then $\|t v\|=0$ for all $t\in \ZZ$. So it is natural to consider non-integer sets of speeds and define
$$
\Rab(n)=\inf_{V\subset \RR\setminus \ZZ\atop |V|=n} \xi(V).
$$
Recall that $\xi(n)$ was defined in a similar way in Section~\ref{subsec:diophantine}, but there the infimum was taken over $V\subset \RR\setminus \QQ$, so 
$$
\Rab(n)\leq \xi(n)=\kappa(n).
$$
Wills~\cite{W1968b} determined $\Rab(n)$ for $n\leq 3$ and Cusick~\cite{C1972} for $n\leq 7$, conjecturing that
$$
\Rab(n)= \frac{1}{w(n)},\quad \text{where }w(n)=\max\{z\in \NN:\, \tfrac{1}{2}\varphi(z)+h(z)\leq n\},
$$
where $\varphi(z)$ is Euler's totient function, $h(z)=0$ if $z$ is prime, and $h(z)$ equals the number of distinct prime divisors of $z$ if $z$ is composite. The conjecture was proved shortly later by Schark~\cite{S1974}. For tight instances, see~\cite{BS2024,C1972}.
Asymptotic estimates are given in~\cite{SW1973},
$$
\Rab(n)\sim \frac{e^{-2\gamma}}{n\log\log n},
$$
where $\gamma= 0.57721$ is the Euler-Mascheroni constant. In contrast to the \LRC, the gap of loneliness in the Rabbit problem turns out to be much narrower.


\subsection{Function fields}\label{sec:ff}

Following the track of establishing $q$--analogs, or function field analogs, of problems in the integers or in the real numbers, Chow and Rimanic~\cite{CR2019} have posed the following function field analog of the \LRC. 

Let $q$ be a prime power and consider the ring of polynomials $\FF_q[X]$ with coefficients in the finite field $\FF_q$ (which plays the role of integers in this analogy), and the field extension $\FF_q((X))$ of Laurent series with coefficients in $\FF_q$ (which plays the role of $\RR$). An element $\alpha=\sum_{i=-\infty}^n \alpha_iX^i\in \FF_q((X))$ can be written as
$$
\alpha =[\alpha]+\|\alpha\|,
$$
where $[\alpha]\in \FF_q[X]$ is a polynomial and $\|\alpha\|=\sum_{i<0}\alpha_iX^i$, the function field analogs of the integer and fractional parts of a real number.  The fractional part belongs to the analog of the {interval $(0,1)$}
, which they denote as
$$
\TT=\{\alpha \in \FF_q((X)): \text{ord} (\alpha){<} 1\},
$$
where $\text{ord}(\alpha)$ is the greatest integer $i\leq n$ such that $\alpha_i\neq 0$.
The distance to the origin of $\|\alpha\|$ is then defined to be $|\alpha|=q^{\text{ord} (\|\alpha\|)}$. For a family $F\subset \FF_q[X]$ denote by
$$
\kappa_q (F)=\sup_{\alpha\in \TT}\min_{f\in F}|\alpha f|.
$$
With this terminology the proposed function field analog of the conjecture reads as follows.

\begin{conjecture}[{Function field} \LRC]\label{conj:ff} For every prime power $q$ and every $F\subset \FF_q[X]\setminus \{0\}$ set of polynomials with 
$$
1\le |F|<\frac{q^{k+1}-1}{q-1} {\quad \text{for some }k\in \NN,}
$$
we have that
$$
\kappa_q (F)\ge q^{-k}.
$$
\end{conjecture}

{The upper bound on the size of $F$ is a necessary condition. Namely, if 
\begin{equation}\label{eq:SPEN}
F=\bigcup_{j=0}^k\{X^{j}+\alpha_{j-1}X^{j-1}+\cdots +\alpha_1X+\alpha_0: \alpha_{0},\alpha_1,\ldots ,\alpha_{j-1}\in \FF_q\},
\end{equation}
then $|F|=\tfrac{q^{k+1}-1}{q-1}$ and $\kappa (F)\le q^{-(k+1)}$.} The authors prove the bound stated in the conjecture for sets $F$ with $|F|\le q^k$. They also prove the conjecture for families $F$ with bounded degree (in analogy of the bounded speeds result in Theorem~\ref{thm:bounded2}) and an approximate version for the case $k=2$.

\section{Collection of Open Problems}\label{sec:OP}

In this section we collect the open problems stated in this survey and propose several new ones, in order to gain further insight into the \LRC. The problems are listed according to their order of appearance in the survey.

\begin{itemize}[leftmargin=0pt]
  \item[] \textbf{Problem 1}: The question of characterizing tight instances of the \LRC\ has seen no further progress since the work of Goddyn and Wong~\cite{GW2006}. The existence of other infinite families of tight instances, besides the ones obtained by appropriately accelerating runners, remains wide open. This question is unsolved even for the case $n=6$, since the proof of the \LRC\ for that case does not provide insight into tight instances.

  \item[] \textbf{Problem 2}: Following the line of work in Section~\ref{sec:gap}, we ask if for every $\epsilon>0$, there exists $n_0\in \NN$ such that for all $n\geq n_0$,
  $$
  \kappa(n)\geq \frac{1+\epsilon}{2n}.
  $$

  \item[] \textbf{Problem 3}: The conjecture is known to hold for any $n\leq 6$. While the proof for $n=6$ is a technical \emph{tour de force}, it is theoretically possible to extend it to other values of $n$ for which $n+1$ is prime; the most natural candidate being $n=10$. In light of recent progress in computer-assisted proofs, it is worth investigating whether such an extension might now be within reach.

  \item[] \textbf{Problem 4}: It is not know if the set of speed instances $V$ for which the conjecture holds is closed under positive translations. Even more, is it true that for any $a\in \NN$ and $n$-set $V$ of positive integers, $\kappa(a+V)\geq \kappa(V)$? The answer is positive if $a$ is sufficiently large with respect to $v_n$, see~\eqref{eq:VNSI}.

  \item[] \textbf{Problem 5}: The reduction by Malikiosis, Santos, and Schymura to a finite set of instances made it possible to prove the conjecture (and its shifted variant) for small values of $n$~\cite{ACS2025,MSS2025}. Pushing the upper bounds on $v_n$ for which the conjecture is known to hold even further remains one of the most promising directions for establishing the conjecture for $n \geq 7$.

  \item[] \textbf{Problem 6}: Tao asked if the \LRC\ is true provided that $v_n\leq 2n$~\cite{T2018}. As discussed in Section~\ref{subsec:lin_bounded}, this is asymptotically true (see Theorem~\ref{thm:bounded3}). On the other hand, Pandey showed that the conjecture holds if $v_n\leq 2n-3$, except for the particular case where $v_1=1$ and $V$ is only contained in an arithmetic progression of length $1$.

  \item[] \textbf{Problem 7}: The Lovász Local Lemma has been used to prove the conjecture given that the sequence is sufficiently lacunary. Variants of the lemma give only limited improvements and it is still open whether there exists $\omega(n)= o(\log{n})$ for which the conjecture holds provided that
  $$
  v_{i+1}\ge \left(1+\frac{\omega(n)}{n}\right)v_i \quad \text{for any }i\in [n-1].
  $$

  \item[] \textbf{Problem 8}: Rifford conjectured (Conjecture~\ref{con:LRC_time}) that if the \LRC\ holds, then there is $N=N(n)$ such that the loneliness of the origin is reached within time $N/v_1$, that is, before the slowest runner has completed $N$ laps around the track. It is a natural idea to use one of the results that reduce the conjecture to a finite number of instances to prove Rifford's conjecture.

  For instance, focus on the approach of Tao~\cite{T2018}. The proof uses induction on $n$ and at each step a runner is added. Let $B(v_1,\dots,v_n;\delta)$ be the $n$-dimensional \emph{Bohr set} of width $\delta>0$, that is
  $$
  B(v_1,\dots,v_n;\delta)= \bigcap_{i=1}^n \{t\in [0,1):\|tv_i\|\leq \delta\}.
  $$
  For every $n\in\NN$, does there exist $N=N(n)$ such that for every $v_1,\dots,v_n$
  \begin{equation}\label{eq:ORJF}
  B\big(v_1,\dots,v_n;\tfrac{1}{n+1}\big) \cap \left[\frac{-N}{v_1},\frac{N}{v_1}\right],
  \end{equation}
  has density bounded below by a positive function depending only on $n$?
  Establishing this would open a path toward proving Rifford's conjecture.

  \item[] \textbf{Problem 9}: In recent years, one of the most actively studied topics related to the \LRP\ has been Kravitz’s Loneliness Spectrum Conjecture. Its amended form (Conjecture~\ref{conj:spec}) seeks to characterize the possible values of the loneliness gap near the tight instances.

  \item[] \textbf{Problem 10}: A paradigm often used in combinatorics is that any object can be partitioned into a \emph{structured} and a \emph{quasirandom} part; see e.g.~\cite{T2008}. Inspired by the results in Section~\ref{sec:random}, we ask whether Theorem~\ref{thm:random} holds for quasirandom sets in the Chung–Graham sense, see e.g.~\cite{CG1992}. In particular, does it hold for Sidon sets, sets including no non-trivial solution of the equation $x+y=z+t$, which are known to be quasirandom?
\end{itemize}

\section{Final comments}\label{sec:final}

In this survey, we aimed to provide an updated and comprehensive overview of the \LRC, hoping to compile and find connections within the majority of the existing literature on this intriguing problem and to inspire further efforts toward its resolution.

Besides the results on small instances (see Section~\ref{sec:small}) of the problem and certain general reductions, the overall picture suggests that the missing part of the proof of the conjecture lies  
\begin{itemize}
    \item[-] between highly structured sets (see Section~\ref{sec:param}) and random sets (see Section~\ref{sec:random});
     \item[-] between compressed sets (see Section~\ref{sec:simple_rem}) and fast growing sets (see Section~\ref{sec:lacunary});
     \item[-] between linearly bounded sets (see Section~\ref{subsec:lin_bounded}) and sets containing a superexponential speed (see Section~\ref{sec:bounded1}).
\end{itemize}
This leaves a vast area of uncertainty. From this perspective, it seems fair to say that the conjecture remains wide open. We hope that this article contributes in keeping the right track towards the solution of this \emph{lovely problem}.

\section*{Acknowledgements} The authors would like to thanks Tom Bohman, Sebastian Czerwiński, Francisco Santos, Matthias Schymura, Jörg M. Wills for their insightful comments that helped improving this survey.  We also thank the anonymous referees for their suggestions.



\begin{thebibliography}{99}
{\small


\bibitem{ACS2025} D. Alcántara, F. Criado, F. Santos, Covering radii of 3-zonotopes and the shifted Lonely Runner Conjecture. arXiv:2506.13379 (2025).

\bibitem{A2013} N. Alon, The chromatic number of random {C}ayley graphs, European J. Combin. 34 (2013) 1232--1243.

\bibitem{AR1999} N. Alon, I.Z. Ruzsa, Non-averaging subsets and non-vanishing transversals, J. Comb. Theory, Ser. A 86 (1) (1999) 1--13.


\bibitem{BS2008} J. Barajas, O. Serra, The lonely runner with seven runners, {Electron. J. Combin.} 15 (2008) R48, 18.

\bibitem{BS2008b} J. Barajas, O. Serra, Distance graphs with maximum chromatic number, {Discrete Math.} 308 (8) (2009) 1355--1365.

\bibitem{BS2009}  J. Barajas, O. Serra, On the chromatic number of circulant graphs, {Discrete Math.} 309 (2008) 5687--5696.

\bibitem{BHS2019} M. Beck, S. Ho\c sten, M. Schymura, Lonely runner polyhedra, {Integers} 19 (2019) A29, 13.

\bibitem{BS2024} M. Beck, M. Schymura, Deep lattice points in zonotopes, lonely runners, and lonely rabbits, {Int. Math. Res. Not. IMRN} 8 (2024) 6553--6578. 

\bibitem{B2017} A. Bernshteyn, The local cut lemma, European J. Combin.63 (2017) 95--114.


\bibitem{BW1972} U. Betke, J.M.Wills, Untere Sehranken fiir zwei diophantische Appro\-xima\-tions-Funktionen, {Monatsh. Math.} 76 (1972) 214--217.

\bibitem{BNV2022}  A. Bhardwaj, V. Narayanan, H. Venkataraman, A few more Lonely Runners, arXiv:2211.08749 (2022).

\bibitem{BGGST1998} W. Bienia, L. Goddyn, P. Gvozdjak, A. Seb\H o, M. Tarsi, Flows, view obstructions, and the lonely runner, {J. Combin. Theory Ser. B} 72 (1998) 1--9.

\bibitem{BHK2001}  T. Bohman, R. Holzman, D. Kleitman, Six lonely runners, {Electron. J. Combin.} 8 (2001) RP3, 49.

\bibitem{BP2022} T. Bohman, F. Peng, Coprime mappings and lonely runners, {Mathematika} 68 (2022) 784--804.



\bibitem{CLZ1999} G.J. Chang, D. Liu, X. Zhu, Distance graphs and $T$-coloring, {J. Combin. Theory (B)} 75 (1999) 159--169.

\bibitem{C1991} Y.G. Chen, On a conjecture in Diophantine approximations II, {J. Number Theory} 37 (2) (1991) 181--198.

\bibitem{C1994} Y.G. Chen, View-obstruction problems in {$n$}-dimensional {E}uclidean space and a generalization of them, {Acta Math. Sinica} 37 (1994) 551--562.

\bibitem{C2000} Y.G. Chen, The best quantitative Kronecker's theorem. {J. London Math. Soc.} 61 (2000) 691--705.

\bibitem{CC1999} Y.G. Chen, T.W. Cusick,  The view-obstruction problem for {$n$}-dimensional cubes, {J. Number Theory} 74 (1999) 126--133.


\bibitem{CG1992} F.R.K. Chung and R.L. Graham, Quasi-random subsets of $\mathbb{Z}_n$, {J. Combin. Theory, Ser. A} 61 (1992) 64--86.

\bibitem{CR2019} S. Chow, L. Rimani\'c, Lonely runners in function fields, {Mathematika} 65 (2019) 677--701.




\bibitem{CMMS2022} J. Cslovjecsek, R.D. Malikiosis, M. Nasz\'odi, M. Schymura, Computing the covering radius of a polytope with an application to lonely runners, {Combinatorica} 42 (2022) 463--490.
 
\bibitem{C1972} T.W. Cusick, Simultaneous diophantine approximation of rational numbers, Acta Arith. 22 (1972) 1--9.

\bibitem{C1973} T.W. Cusick, View-obstruction problems, Aequationes Math. 9 (1973) 165–170.

\bibitem{C1974} T.W. Cusick, View-obstruction problems in $n$-dimensional geometry. {J. Combin. Theory Ser. A} 16 (1974) 1--11. 

\bibitem{C1982} T.W. Cusick, View-obstruction problems {II},  {Proc. Amer. Math. Soc.} 84 (1982) 25--28.

\bibitem{CP1984} T.W. Cusick, C. Pomerance, View-obstruction problems {III},  {J. Number Theory} 19 (1984) 131--139.

\bibitem{C2012} S. Czerwi\'nski, Random runners are very lonely, {J. Combin. Theory Ser. A} 119 (2012) 1194--1199.

\bibitem{C2018} S. Czerwi\'{n}ski, The lonely runner problem for lacunary sequences, {Discrete Math.} 341 (2018) 1301--1306.

\bibitem{CG2008} S. Czerwi\'{n}ski, J. Grytczuk, Invisible runners in finite fields, {Inform. Process. Lett.} 108 (2008) 64--67.


\bibitem{D2011} A. Dubickas, The lonely runner problem for many runners, {Glas. Mat. Ser. III} 46 (2011) 25--30.

\bibitem{DT2020} A. Dumitrescu, C.D. T\'{o}th, Problems on track runners {Comput. Geom.} 88 (2020) 101611.


 \bibitem{EES1985} R.B. Eggleton, P. Erd\H{o}s, D.K. Skilton, Colouring the real line, {J. Combin. Theory (B)} 39 (1985) 86--100.

\bibitem{E1962} P. Erd\H os, On the integers relatively prime to $n$ and on a number-theoretic function considered by Jacobsthal, Math. Scand. 10 (1962) 163--170.

\bibitem{E1975} P. Erd\H{o}s, Problems and results on Diophantine approximations II, Actes Colloq. Marseille-Luminy 1974, Lecture Notes in Math. 475 (1975) 89--99.

\bibitem{FS2023} H.T. Fan, A. Sun, Amending the Lonely Runner Spectrum Conjecture, arXiv:2306.10417 (2023).


\bibitem{FF1958} L.R. Ford, D. R. Fulkerson, Network flow and systems of representatives, Canadian Journal of Mathematics 10 (1958) 78-84.


\bibitem{GK2025} V. Giri, N. Kravitz,  The structure of Lonely Runner spectra, arXiv:2304.01462 (2025).



\bibitem{GW2006} L. Goddyn, E.B. Wong, Tight instances of the lonely runner, {Integers} 6 (2006) A38, 14.

\bibitem{GR2022} F. Gon\c{c}alves, J.P.G. Ramos, Bounds for the lonely runner problems via linear programming, {Bull. Braz. Math. Soc. (N.S.)} 53 (2022) 595--603.

\bibitem{G2017} B. Green, On the chromatic number of random {C}ayley graphs,  {Combin. Probab. Comput.}  26 (2017) 248--266.

\bibitem{GK1987} P.M. Gruber, C.G. Lekkerkerker, \emph{Geometry of Numbers} 2 edn., North-Holland Mathematical Library 37 (1987).

\bibitem{H1976b} N.M. Haralambis, Sets of integers with missing differences, {J. Combin. Theory Ser. A} 23 (1) (1977) 22--33.

\bibitem{HW1979} G.H. Hardy,  E.M. Wright, {An introduction to the theory of numbers}, Oxford university press 1979.

\bibitem{HM2017} M. Henze, R.D. Malikiosis, On the covering radius of lattice zonotopes and its relation to view-obstructions and the lonely runner conjecture, {Aequationes Math.} 91 (2017) 331--352.


\bibitem{H1976} D. Hunter,  An Upper Bound for the Probability of a Union, Journal of Applied Probability 13 (3) (1976) 597--603.

\bibitem{J1961} E. Jacobsthal, \"Uber {S}equenzen ganzer {Z}ahlen, von denen keine zu {$n$} teilerfremd ist. {I}, {II}, {III}, Norske Vid. Selsk. Forh. 33 (1961) 117--139.

\bibitem{K2001} Y. Katznelson, Chromatic numbers of Cayley graphs on $\ZZ$ and recurrence, Combinatorica 21 (2001) 211--219.

\bibitem{JK2024} V. Jain, N. Kravitz, Relative Lonely Runner spectra, arXiv:2411.12684 (2024).



\bibitem{K2021} N. Kravitz, Barely lonely runners and very lonely runners: a refined approach to the lonely runner problem,  {Comb. Theory} 1 (2021) P17, 24.

\bibitem{L2008} D. Liu, From rainbow to the lonely runner: a survey on coloring parameters of distance graphs, {Taiwanese J. Math.} 12 (2008) 851--871.

\bibitem{LR2020} D. Liu, G. Robinson, Sequences of integers with three missing separations, {European J. Combin.} 85 (2020) 103056.

\bibitem{LS2013} D. Liu, A. Sutedja, Chromatic number of distance graphs generated by the sets {$\{2,3,x,y\}$}, {J. Comb. Optim.} 25 (2013) 680--693.

\bibitem{LZ2004} D. Liu, X. Zhu, Fractional chromatic number and circular chromatic number for distance graphs with large clique size, {J. Graph Theory} 47 (2004) 129--146.

\bibitem{M1980} B. de Mathan, Numbers contravening a condition in density modulo 1, Acta Math. Hungar. 36 (1980) 237--241.

\bibitem{MSS2025} R. D. Malikiosis, F. Santos, M. Schymura
Linearly-exponential checking is enough for the Lonely Runner Conjecture and some of its variants,  arXiv:2411.06903v1 [math.CO] 

\bibitem{P2009} R.K. Pandey, A note on the lonely runner conjecture, {J. Integer Seq.} 12 (2009) A9 4.

\bibitem{P2010} R.K. Pandey, On the lonely runner conjecture, {Math. Bohem.} 135 (2010) 63--68.

\bibitem{PR2023} R.K. Pandey, N. Rai, Maximal density and the kappa values for the families $\{a,a+1,2a+1,n\}$ and $\{a,a+1,2a+1,3a+1,n\}$, {Math. Slovaca} 73 (3) (2023) 643--656.

\bibitem{PS2022} R.K. Pandey, A. Srivastava, Maximal density of integral sets with missing differences and the kappa values, {Taiwanese J. Math.} 26 (2022) 17--32.

\bibitem{PT2011} R.K. Pandey, A. Tripathi, On the density of integral sets with missing differences from sets related to arithmetic progressions, {J. Number Theory} 131 (2011) 634--647.

\bibitem{PS2016} G. Perarnau, O. Serra, Correlation among runners and some results on the lonely runner conjecture, {Electron. J. Combin.} 23 (2016) P50, 22.

\bibitem{PS2010} Y. Peres, W. Schlag.  Two Erdős problems on lacunary sequences: chromatic number and Diophantine approximation. {Bull. Lond. Math. Soc.} 42 (2) (2010) 295--300.

\bibitem{P1979} A.D. Pollington, On the density of the sequence $\{n_k\xi\}$, Illinois J. Math. 23 (1979) 511--515.

\bibitem{P2022} C. Pomerance, Coprime matchings,  {Integers} 22 (2022) A2 9.

\bibitem{PS1980} C. Pomerance, J.L Selfridge, Proof of DJ Newman’s coprime mapping conjecture. {Mathematika} 27 (1) (1980) 69--83.

\bibitem{R2004} J. Renault, View-obstruction: a shorter proof for 6 lonely runners, {Discrete Math.} 287 (2004) 93--101.

\bibitem{R2022} L. Rifford, On the time for a runner to get lonely, {Acta Appl. Math.} 180 (2022) P15 66.

\bibitem{RTV2002} I.Z. Ruzsa, Zs. Tuza, M. Voigt, Distance graphs with finite chromatic number, {J. Combin. Theory Ser. B} 85 (1) (2002) 181--187.

\bibitem{S1974} R. Schark, Eine diophantische Approximations-Funktion, Monatsh. Math. 78 (1974) 131--146.

\bibitem{SW1973} R. Schark, J.M. Wills, Asymptotisches Verhalten einer diophantischen Approximationsfunktion, Acta Arith. 22 (1973) 129--136.

\bibitem{S1976} I.J. Schoenberg, Extremum problems for the motions of a billiard ball II, The $L_{\infty}$ norm, Nederl. Akad. Wetensch. Proc. Ser. A 79 Indag. Math. 38 (3) (1976), 263--279.

\bibitem{SW2018} M. Schymura, J.M. Wills, Der einsame Läufer, \emph{Mitteilungen der Deutschen Mathematiker-Vereinigung} 26 (1) (2018) 14--17.

\bibitem{S1981} P.D. Seymour, Nowhere-zero 6-flows. {J. Combin. Theory Ser. B} 30 (2) (1981) 130--135.



\bibitem{T2008} T. Tao, Structure and Randomness, American Mathematical Society (2008).

\bibitem{T2018} T. Tao, Some remarks on the lonely runner conjecture, {Contrib. Discrete Math.} 13 (2018) 1--31.

\bibitem{W1916} H. Weyl, Uber die Gleichverteilung von Zahlen modulo Eins, {Mathematische Annalen} 77 (1916) 313--352. 


\bibitem{W1965} J.M. Wills, Zwei Probleme der inhomogenen diophantischen Approximation. PhD thesis (1965), TU Berlin.

\bibitem{W1967} J.M. Wills, Zwei Sätze über inhomogene diophantische Approximation von Irrationalzahlen. {Monatsh. Math.} 71 (1967) 263--269.


\bibitem{W1968} J.M. Wills, Zur simultanen homogenen diophantischen Approximation. I, {Monatsh. Math.} 72 (1968) 254--263.

\bibitem{W1968b} J.M. Wills, Zur simultanen homogenen diophantischen Approximation. II, {Monatsh. Math.} 72 (1968) 368--381.

\bibitem{WL2004} J. Wu, W. Lin, Circular chromatic numbers and fractional chromatic numbers of distance graphs with distance sets missing an interval, {Ars Combin.} 70 (2004) 161-168.

\bibitem{YPW2020} Q.H.Yang, T. Pan, J.D. Wu, On optimal {$M$}-sets related to {M}otzkin's problem, {J. Math.} (2020) 7457625 5.

\bibitem{Z2001} X. Zhu, Circular chromatic number: a survey. Combinatorics, graph theory, algorithms and applications {Discrete Math.} 229 (1--3) (2001) 371--410.



}

\end{thebibliography}
\end{document}